\documentclass{amsart}
\usepackage[utf8]{inputenc}
\usepackage{fancyhdr}
\usepackage{chngcntr}
\usepackage{mathabx}
\usepackage{mathtools}

\usepackage{enumitem}
\usepackage{breakurl}
\usepackage[colorlinks, breaklinks, linkcolor=blue]{hyperref}

\usepackage[ngerman,american]{babel}
\usepackage[T1]{fontenc}
\usepackage{upgreek}
\usepackage{dsfont}
\usepackage{suffix}
\newcommand\be[4][{}]{\ifthenelse{\equal{#1}{}}{\begin{#2}\label{#3}#4\end{#2}}{\begin{#2}[#1]\label{#3}#4\end{#2}}}
\WithSuffix\newcommand\be*[3][{}]{\ifthenelse{\equal{#1}{}}{\begin{#2}#3\end{#2}}{\begin{#2}[#1]#3\end{#2}}}

\usepackage{xstring}

\def \LEMMA {l}

\newcommand{\labelmap}[1]{
	\IfEqCase{#1}{
		{lemma}{\LEMMA}
		{equation}{eq}
		{corollary}{c}
		{proposition}{p}
		{theorem}{t}
		{remark}{r}
		{example}{ex}
		{section}{s}
		{figure}{f}
	}[\PackageError{labelmap}{Undefined option to labelmap: #1}{}]
}

\newcommand{\beltest}[4][{}]
{
	\ifthenelse{\equal{#1}{}}
	{\begin{#2}\label{\labelmap{#2}:#3}#4\end{#2}}
	{\begin{#2}[#1]\label{\labelmap{#2}}#4\end{#2}}
}

\newcommand{\ft}[1]{\text{#1}} 
\newcommand{\ftqq}[1][and]{\quad\ft{#1}\quad}

\newif\ifSIMPLE

\usepackage{xspace}
\newcommand{\iww}[1][take]{in which we #1\xspace}
\newcommand{\lhs}{left-hand side\xspace}
\newcommand{\rhs}{right-hand side\xspace}
\newcommand{\wrt}{with respect to\xspace}

\newcommand{\mb}[2][]{\ifthenelse{\equal{#1}{}}{\left(#2\right)}{\big(#2\big)}} 

\WithSuffix\newcommand\mb*[1]{(#1)} 
\newcommand{\mbs}[1]{\left\{#1\right\}} 
\newcommand{\mblo}[1]{\left(#1\right]} 
\newcommand{\mbb}[1]{\left[#1\right]} 
\newcommand{\mbm}[1]{\left|#1\right|} 
\newcommand{\mbfl}[1]{\left\lfloor#1\right\rfloor}

\newcommand{\mf}[2]{\frac{#1}{#2}} 

\newcommand{\on}[1]{\operatorname{#1}}

\newcommand{\iunit}{\on{i}}

\newcommand{\SR}[1][{ }]{\ft{#1}|\ft{#1}}

\usepackage{amsthm}

\newcommand{\ch}[2][]{\ifthenelse{\equal{#1}{}}{Chapter \ref{ch:#2}}{Chapter \ref{ch:#2:#1}}}
\newcommand{\se}[2][]{\ifthenelse{\equal{#1}{}}{Section \ref{s:#2}}{Section \ref{s:#2:#1}}}

\newcommand{\rem}[2][]{\ifthenelse{\equal{#1}{}}{Remark \ref{r:#2}}{\tm{#2} \ref{r:#2:#1}}}

\newcommand{\tm}[2][]{\ifthenelse{\equal{#1}{}}{Theorem \ref{t:#2}}{\tm{#2} \ref{t:#2:#1}}}
\WithSuffix\newcommand\tm*[2][]{\ifthenelse{\equal{#1}{}}{\ref{t:#2}}{\ref{t:#2:#1}}}

\newcommand{\lem}[2][]{\ifthenelse{\equal{#1}{}}{Lemma \ref{l:#2}}{\lem{#2} \ref{l:#2:#1}}}
\WithSuffix\newcommand\lem*[2][]{\ifthenelse{\equal{#1}{}}{\ref{l:#2}}{\ref{l:#2:#1}}}

\newcommand\co[2][]{\ifthenelse{\equal{#1}{}}{Corollary \ref{c:#2}}{\co{#2} \ref{c:#2:#1}}}
\WithSuffix\newcommand\co*[2][]{\ifthenelse{\equal{#1}{}}{\ref{c:#2}}{\ref{c:#2:#1}}}

\newcommand{\eq}[2][]{\ifthenelse{\equal{#1}{}}{Equation \ref{eq:#2}}{\eq{#2} \ref{eq:#2:#1}}}
\WithSuffix\newcommand\eq*[2][]{\ifthenelse{\equal{#1}{}}{(\ref{eq:#2})}{(\ref{eq:#2:#1})}}

\newcommand{\qe}[2][]{\ifthenelse{\equal{#1}{}}{Question \ref{qe:#2}}{\eq{#2} \ref{qe:#2:#1}}}
\WithSuffix\newcommand\qe*[2][]{\ifthenelse{\equal{#1}{}}{(\ref{eq:#2})}{(\ref{qe:#2:#1})}}

\newcommand{\pb}[2][]{\ifthenelse{\equal{#1}{}}{Problem \ref{pb:#2}}{\eq{#2} \ref{pb:#2:#1}}}
\WithSuffix\newcommand\pb*[2][]{\ifthenelse{\equal{#1}{}}{(\ref{pb:#2})}{(\ref{pb:#2:#1})}}

\newcommand{\ex}[2][]{\ifthenelse{\equal{#1}{}}{Example \ref{ex:#2}}{\ex{#2} \ref{ex:#2:#1}}}

\newcommand{\de}[2][]{\ifthenelse{\equal{#1}{}}{Definition \ref{d:#2}}{\de{#2} \ref{d:#2:#1}}}

\newcommand{\cnd}[2][]{\ifthenelse{\equal{#1}{}}{Condition \ref{cnd:#2}}{\de{#2} \ref{cnd:#2:#1}}}

\newcommand{\pp}[2][]{\ifthenelse{\equal{#1}{}}{Proposition \ref{p:#2}}{\pp{#2} \ref{p:#2:#1}}}
\WithSuffix\newcommand\pp*[2][]{\ifthenelse{\equal{#1}{}}{\ref{p:#2}}{\ref{p:#2:#1}}}

\newlist{enumroman}{enumerate}{1}
\setlist[enumroman,1]{label=(\roman*),before*=\singlespacing}


\newcommand\eqn[2][]{\ifthenelse{\equal{#1}{}}{(\ref{eq:#2})}{(\ref{eq:#2:#1})}}
\newcommand{\ppn}[2][]{\ifthenelse{\equal{#1}{}}{\ref{p:#2}}{\ref{p:#2:#1}}}

\newcommand{\MB}[1]{\mathbb{#1}}

\usepackage[backend=biber, isbn=false]{biblatex}
\renewbibmacro{in:}{}
\addbibresource{../../../bibtex/math/math.bib}
\AtEveryBibitem
{
	\clearlist{address}
	\clearlist{language}
	\clearfield{month}
}
\usepackage{abstract}
\setlist[enumerate,1]{label=(\roman*),font=\normalfont}

\setcounter{section}{0}

\setlength\parindent{0pt}

\newcommand{\Addresses}{{
		\bigskip
		\footnotesize
		
		Alexander Adam: \textsc{Max Planck Institute for Mathematics, Vivatsgasse 7, 53113 Bonn, Germany}\par\nopagebreak
		\textit{E-mail:} \texttt{a.adam.sci@protonmail.com}
	}}

\be*{document}
{
	\title[Hypergeometric expansion related to the Hurwitz zeta function]{Generalized hypergeometric expansion related to the Hurwitz zeta function}
	\author{Alexander Adam}
	\address{\Addresses}
	\date{\today}

	\maketitle
	\be*{abstract}
	{
		We study the incomplete Mellin transformation of the fractional part and the related log-sine function when composed by an affine complex map. We evaluate the corresponding integral in two different ways which yields equalities with series in hypergeometric functions on each side. These equalities capture basic analytic properties of the Hurwitz zeta function like meromorphic extension and the functional equation. We give rates of convergence for the involved series. As a special case we find that the exponential rate of convergence for the deformation of the log-sine to the fractional part in the imaginary component is carried over to these equalities.
	}

	
	\noindent
	\tableofcontents
	\section{Introduction}
	The Riemann zeta function
	\be{align}{eq:zeta}
	{
		\zeta\mb{s}\coloneqq\sum_{k=1}^\infty k^{-s},\quad s\in\MB{C} \ft{ such that } \Re s>1,
	}
	is arguably one of the best studied Dirichlet series. We can easily rewrite
	\[\pi\mf{\zeta\mb{s}}{\sin\mb{\pi s}}\]
	(or any similar Dirichlet series) as a Gauss hypergeometric expansion, utilizing the remarkable identity
	\be{align}{eq:hyperidentity}
	{
		\mf\pi{\sin\mb{\pi s}}z^s=z\mf{{}_2F_1\mb{1,1-s;2-s,-z}}{1-s}+\mf{{}_2F_1\mb{1,s;1+s,-\mf 1z}}{s},
	}
	which arises as a special case as one of the $20$ linear relations of Kummmer's 24 solutions of the hypergeometric differential equation \cite[Section 3.8-3.9]{Luke_1969}, where ${}_2F_1$ denotes the Gauss hypergeometric function. Notably the above identity is valid for every $z\in\MB{C}\setminus\mbb{-1,0}$ and every $s\in\MB{C}\setminus\MB{Z}$ if the branch for the logarithm is taken the same on both sides in \eq*{hyperidentity}.\\
	That way $\zeta$ is defined in \eq*{zeta} the convergence is only on the open half plane $\Re s>1$ and the na\"ive use of the identity in \eq*{hyperidentity} is of no further help. A typical method for extending $\zeta$ meromorphically on $\MB{C}$ is to derive an extension on $\Re s>0$ and to derive an analytic expression relating $\zeta\mb{s}$ with $\zeta\mb{1-s}$ on the open strip $0<\Re s<1$. The existence of a meromorphic continuation for $\zeta$ on $\MB{C}$ and the evaluation at special arguments for $\zeta$ has led to the writing of e.g.
	\[1+2+3+\dots=-\frac 12=\zeta\mb{-1},\]
	as a short-hand which is of course wrong in the strict sense as the limit $\lim_{s\to -1}\zeta\mb{s}$ cannot be interchanged with the sum since the Dirichlet series is not defined at $-1$, although its meromorphic continuation is.\\
	We present here a simple pattern into which falls the meromorphic continuation of $\zeta$ and where an exchange of the limit in $s$ is possible. This simple pattern consists of a family of functions $f_k\mb{s}$ and $g_k\mb{s}$, where $k\in\MB{N}$ and $s\in \MB{C}$ and fixes the analytic properties of the sum
	\[\sum_{k=1}^\infty f_k\mb{s}+g_k\mb{s},\]
	depending on $s$. For $\Re s<0$ it is required that the sums
	\[\sum_{k=1}^\infty f_k\mb{s}\ftqq \sum_{k=1}^\infty g_k\mb{s}\]
	to converge absolutely to holomorphic functions in $s$ and uniform in $s$ on certain domains in $\MB{C}$. Instead, for $s\in\MB{C}$ the same kind of convergence is only required to hold for the sum
	\[\sum_{k=1}^\infty f_k\mb{s}+g_k\mb{s}.\]
	Hence the real part of $s$ serves as an interpolation quantity for the separability of the terms $f_k\mb{s}+g_k\mb{s}$, $k\in\MB{N}$, \wrt absolute convergence.\\	
	A prototypical example poses the Dirichlet eta function for every $s\in\MB{C}$, where $\Re s>0$,
	\[\eta\mb{s}\coloneqq-\sum_{k=1}^\infty\mb{-1}^kk^{-s}.\]
	If we set for every $k\in\MB{N}$
	\[f_k\mb{s}\coloneqq \mb{2k-1}^{-s}\ftqq g_k\mb{s}\coloneqq -\mb{2k}^{-s}\]
	then we can write
	\[\eta\mb{s}\coloneqq\sum_{k=1}^\infty f_k\mb{s}+g_k\mb{s}.\]
	Here $\sum_{k=1}^\infty \mbm{f_k\mb{s}}$ and $\sum_{k=1}^\infty \mbm{g_k\mb{s}}$ exist for every $\Re s>1$ and $\sum_{k=1}^\infty \mbm{f_k\mb{s}+g_k\mb{s}}$ exists for every $0<\Re s<\infty$. For $\Re s>1$ the termwise convergence allows then to relate $\eta$ to $\zeta$ via the functional equation
	\[\eta\mb{s}=\sum_{k=1}^\infty f_k\mb{s}-g_k\mb{s}+2g_k\mb{s}=\mb{1-2^{1-s}}\zeta\mb{s}.\]
	We argue that an analogous pattern appears naturally also directly for $\zeta$ as it follows from a straight-forward analysis of a well-known and canonical integral representation for $\zeta$ for every $\Re s>0$
	\be{align}{eq:mellinzeta}
	{
		\zeta\mb{s}-\mf s{s-1}+\mf 12=s\int_1^\infty\mb{\mf 12-\mbs{x}}x^{-s-1}\operatorname{d}x,
	}
	where $\mbs{\cdot}$ denotes the fractional part on $\MB{R}$ (\cite[Section 14.1]{Baker_2012}).
	Therefore we stress that this pattern reveals how the machinery of meromorphic continuation for $\zeta$ works through additive termwise deformation of its defining Dirichlet series. Such a deformation is of course not unique but the deformation resulting from \tm{maintheorem} below follows from a direct calculation. In particular, what underlies here the general convergence behavior is the smoothing behavior of the Mellin transform, noting that the Fourier series for almost all $x\in\MB{R}\setminus\MB{Z}$
	\be{align}{eq:fourierseries}
	{
		\mf 12-\mbs{x}=\mf 1\pi\sum_{k=1}^\infty\mf{\sin\mb{2\pi xk}}k
	}
	is not absolutely convergent. A (very) small search through the literature \cite{Baker_2012, Patterson_1988, Titchmarsh_1986} seems to reveal that the integral
	\[\int_0^1 \sin\mb{x}x^{-s-1}\operatorname{d}x\]
	has never been evaluated there in the context of the functional equation for $\zeta$. While this defines an entire function and gives rise to the hypergeometric function ${}_1F_2$ with defining series given in \eq*{hypergeo12} below, instead, in the derivation of the functional equation the completed integral
	\be{align}{eq:gammaintegral}
	{
		\int_0^\infty \sin\mb{x}x^{-s-1}\operatorname{d}x=-\sin\mb{\mf  {\pi s}2}\Gamma\mb{-s}.
	}
	is directly used. Moreover, the recurrence relation for the values $\zeta\mb{2n}$, $n\in\MB{N}$, from \lem{recurrence:odd} follows (from) a recurrence relation for ${}_1F_2$. The converse result for $\zeta\mb{2n+1}$ from \lem{recurrence:even} follows analogously but the closed evaluation of the \lhs remains unknown but is provided by the result from \lem{mellintransform:III} as a hypergeometric expansion. Regarding this problem, generalized hypergeometric series have been used in order to find (computational) good representations for such values or for the computation of related special functions as a whole \cite{Chu_1997,Bin_2009,Crandall_2012, Borwein_2013, Coffey_2014,Sofo_2015}.\\
	A slight generalization for $\zeta$ poses the Hurwitz zeta function for every $x>0$
	\[\zeta\mb{s,x}\coloneqq\sum_{k=0}^\infty \mb{k+x}^{-s}\ftqq\ft{for every }\Re s>1.\]
	It is known that this function admits a meromorphic continuation in $s\in\MB{C}$ and satisfies the Hurwitz identity.
	For every $\alpha,z\in\MB{C}$ and $\beta,\gamma\in\MB{C}\setminus \MB{Z}_{\le 0}$ the hypergeometric function ${}_1F_2$ is defined as the exponential generating series
	\be{align}{eq:hypergeo12}
	{
		{}_1F_2\mb{\alpha;\beta,\gamma;z}\coloneqq \sum_{n=0}^\infty \mf{\mb{\alpha}_n}{\mb{\beta}_n\mb{\gamma}_n}\mf {z^n}{n!},
	}
	where the the bracket $\mb{\cdot}_n$ is understood as a quotient of two Gamma functions
	\[\mb{\cdot}_n\coloneqq \mf{\Gamma\mb{\cdot+n}}{\Gamma},\]
	also known as the Pochhammer symbol. We are able to recover the meromorphic continuation and the functional equation for the Hurwitz zeta function from a single equation given by the following theorem:
	\be{theorem}{t:maintheorem}
	{
		Let $s\in\MB{C}$ such that $\Re s>1$, let $a\in\MB{R}\setminus\mbs{0}$, $b\in\MB{R}$ and let $d>0$. Setting
		\[R_{a,b,d}\mb{s}\coloneqq s\mb{\int_1^{\mf{1-\mbs{b}}a}-\int_1^d}\mbfl{ax+\mbs{b}}x^{-s-1}\on{d}x,\]
		then it holds
		\be*{align*}
		{	
			\mf \pi s&\be*{cases}
			{
				a^s\zeta\mb{s,1-\mbs{b}}-\mb{\mf {ads}{s-1}+\mbs{b}-\mf 12}d^{-s}+R_{a,b,d}\mb{s},&\ft{if }a>0\\
				-\mb{-a}^s\zeta\mb{s,1-\mbs{-b}}+\mb{\mf {ads}{1-s}+\mbs{-b}-\mf 12}d^{-s}-R_{-a,-b,d}\mb{s},&\ft{if }a<0
			}=\\
			&\sum_{k=1}^\infty\mf{\sin\mb{2\pi bk}}{d^ssk}{}_1F_2\mb{-\mf s2;\mf 12,1-\mf s2;-\pi^2a^2d^2k^2}\\
			&-\mf{2\pi ad^{1-s}\cos\mb{2\pi bk}}{1-s}{}_1F_2\mb{\mf 12\mb{1-s};\mf 32,1+\mf 12\mb{1-s};-\pi^2a^2d^2k^2}\\
			&+k^{s-1}\mb{2\pi}^s\Gamma\mb{-s}\mf 1{2\iunit}\mb{\mb{-\iunit a}^se^{\iunit 2\pi{bk}}-\mb{\iunit a}^se^{-\iunit 2\pi{bk}}},
		}
		where the series on the \rhs converges to a holomorphic function in $s\in\MB{C}$. The convergence is absolute and uniform in $b$, in $a,d$ bounded away from $0$ and $\infty$ and for every $\epsilon>0$ in $s$ with $\Re s>-\epsilon$ and $\mbm{\Im s}<\epsilon$. If $\Re s<1$ and $b\in\MB{R}\setminus\MB{Z}$ or $\Re s<0$ and $b\in\MB{Z}$ it holds
		\be*{align*}
		{
			\sum_{k=1}^\infty&\mf{2\pi a\cos\mb{2\pi bk}}{1-s}{}_1F_2\mb{\mf 12\mb{1-s};\mf 32,1+\mf 12\mb{1-s};-\pi^2a^2k^2}\\
			-&\mf{\sin\mb{2\pi bk}}{sk}{}_1F_2\mb{-\mf s2;\mf 12,1-\mf s2;-\pi^2a^2k^2}=\mf \pi{s}\mb{\mf {as}{s-1}+b-\mf 12+R_{a,b,1}\mb{s}}.
		}
	}
	Note that the above theorem gives immediately the holomorphic continuation in $s~\in~\MB{C}$ for $\zeta\mb{s,1-\mbs{b}}+\mb{1-s}^{-1}$, the functional equation for $\zeta\mb{s,1-\mbs{b}}$, the value $\zeta\mb{0,1-\mbs{b}}=\mbs{b}-\mf 12$ and expressions for
	\[\mb{\partial_s\zeta}\mb{0,1-\mbs{b}}\ftqq
	\mb{\zeta\mb{s,1-\mbs{b}}+\mb{1-s}^{-1}}_{|s=1},\]
	calculating the limits $s\to 0$ and $s\to 1$, respectively, as well. Moreover, utilizing the identity in \eq*{hyperidentity}, we can cast the \rhs completely as a series in sums of (generalized) hypergeometric functions. The special case $b=0$ belongs to the Riemann zeta function. In this case, if we multiply both sides by $\mb{\mb{2\pi}^s\Gamma\mb{-s}\sin\mb{\mf \pi 2s}}^{-1}$, substituting $s\mapsto 1-s$ we recover on the \rhs an additive deformation of $k^{-s}$ and the series shows an analogous convergence behavior as mentioned above for the prototypical example.\\
	The converse result for \tm{maintheorem} arises by considering the exchange $\sin\mapsto \cos$, $\cos\mapsto -\sin$ on the \rhs of the principal identity from \tm{maintheorem}. The \lhs in the identity from the following theorem takes on a more complicated shape and the Gauss hypergeometric function appears: 
	\be{theorem}{t:maintheorem:B}
	{
		Let $s\in\MB{C}$ such that $\Re s>1$, let $a\in\MB{R}\setminus\mbs{0}$, $b\in\MB{R}$, $d>0$ and let $n\in\MB{Z}$ be minimal such that $n>\mbm{a}d+b\on{sgn}{a}$. Setting
		\be*{align*}
		{
			H_{a,b,d}\mb{s}\coloneqq&\mf{\log 2}{sd^s}+1_{|\MB{R}\setminus\MB{Z}}\mb{ad+b}\mf {\log\mb{\sin\mb{\pi\mb{ad+b}}}^2}{2sd^s}+1_{|\MB{Z}}\mb{b}\mf {1}{s^2d^s}\\
			&+\mf {1}{sd^s}1_{|\MB{Z}}\mb{ad+b}\mb{\mf {1}{s}-{\psi\mb{1+s}-\gamma}+\mf 12\log\mb{\pi ad}^2}\\
			&-\mf{1}{d^{s}}\sum_{\substack{k\in\MB{Z}\\
					k\not\in\mbs{b,ad+b}}}\mf{ad}{b-k}\mf{{}_2F_1\mb{1,1-s;2-s;-\mf{ad}{b-k}}}{s\mb{1-s}}\\
			&+\iunit \mf{\pi\mbm{a}^s}{s}\mb{\zeta\mb{s,1-\mbs{b\on{sgn}a}}-\zeta\mb{s,n-b\on{sgn}a}},
		}
		where $\gamma$ denotes the Euler-Mascheroni constant and $\psi$ the Digamma function, then it holds
		\be*{align*}
		{
			H_{a,b,d}\mb{s}&+\mf{\pi\mbm{a}^s}{s}\mb{\mf{\zeta\mb{s,1-\mbs{-b\on{sgn}a}}}{\sin\mb{\pi s}}+\mf{\zeta\mb{s,1-\mbs{b\on{sgn}a}}}{\tan\mb{\pi s}}}=\\
			&-\sum_{k=1}^\infty\mf{\cos\mb{2\pi bk}}{d^ssk}{}_1F_2\mb{-\mf s2;\mf 12,1-\mf s2;-\pi^2a^2d^2k^2}\\
			&+\mf{2\pi ad^{1-s}\sin\mb{2\pi bk}}{1-s}{}_1F_2\mb{\mf 12\mb{1-s};\mf 32,1+\mf 12\mb{1-s};-\pi^2a^2d^2k^2}\\
			&+k^{s-1}\mb{2\pi}^s\Gamma\mb{-s}\mf 1{2}\mb{\mb{-\iunit a}^se^{\iunit 2\pi{bk}}+\mb{\iunit a}^se^{-\iunit 2\pi{bk}}}.
		}
		The series on the \lhs converges to a meromorphic function (with poles in $s\in\MB{N}_{\ge 2}$) and the summation in $k$ is understood from $0$ to $\infty$ in absolute value. The convergence on the \lhs is absolute and uniform in $s$ for every $\epsilon>0$ such that $\Re s <\epsilon$ and $s$ bounded away from $\MB{N}_{\le \epsilon}$, uniform in $b$ uniformly bounded away from $\MB{Z}$ and uniform in $a$ and $d$ with $ad+b$ bounded away from $\MB{Z}$ and $ad$ bounded. The series on the \rhs converges to a holomorphic function in $s$. 
		The convergence is absolute and uniform in $b$, in $a,d$ bounded away from $0$ and $\infty$ and for every $\epsilon>0$ in $s$ with $\Re s>-\epsilon$ and $\mbm{\Im s}<\epsilon$.\\
		If $\Re s<1$ and $b\in\MB{R}\setminus\MB{Z}$ or $\Re s<0$ and $b\in\MB{Z}$ it holds
		\be*{align*}
		{
			\sum_{k=1}^\infty&\mf{2\pi a\cos\mb{2\pi bk}}{1-s}{}_1F_2\mb{\mf 12\mb{1-s};\mf 32,1+\mf 12\mb{1-s};-\pi^2a^2k^2}\\
			+&\mf{\cos\mb{2\pi bk}}{sk}{}_1F_2\mb{-\mf s2;\mf 12,1-\mf s2;-\pi^2a^2k^2}=-H_{a,b,1}\mb{s}.
		}
	}
	As we have already explained \tm{maintheorem}-\tm*{maintheorem:B} appear naturally in the study of the Hurwitz zeta function and both theorems can be superimposed to one equality which \rhs is \tm{convergence} presented in \se{result}. In principle, \tm{maintheorem:B} can be extended to all $a\in\MB{C}\setminus\mbs{0}$, $b\in\MB{C}$ by combining \tm{convergence} with \lem{mellintransform:IV} (but we do not write down the identity explicitly as we have no immediate use in the direction of \co{limit} below). We then discuss the meaning of the integrand in the underlying integral transform, apart from the rather obvious connection to Poisson summation (as this leads to a fused interpretation for integrand and integral kernel), as well as possible extensions of our result.\\ 
	Then \se{proof} is completely devoted to the proofs which is then followed by \se{app} where we present some immediate consequences which are in principal all known but arise here either as direct consequence of the presented theorems or as by-products along their proofs.
	\section{Result}\label{s:result}
	\setcounter{theorem}{1}
	As we explained in the introduction, our basic starting point for \tm{maintheorem} in the most simple case ($a=d=1$, $b=0$) is the integral in \eq*{mellinzeta}.
	The principal identity follows then by evaluation of the integral in two different ways. We proceed analogously for \tm{maintheorem:B}.
	In particular, the evaluation underlying the \rhs in the identities from \tm{maintheorem}-\tm*{maintheorem:B} follow a unique algorithm and the result is captured in \tm{convergence} below.\\		
	We endow the logarithm with the principal branch, meaning for every $t\in\mblo{-1,1}$ it holds
	\[\log e^{\iunit \pi t}=\iunit \pi t,\]
	and for every $a,b\in\MB{C}$, $a\neq 0$, we set
	\be{align}{eq:complexpow}
	{
		a^b\coloneqq e^{b\log a}.
	}
	This is important, regarding the extension of the identity in \eq*{hyperidentity} to all $\MB{C}\setminus\mbs{-1,0}$.
	Moreover, for every $z\in\MB{C}$ we set as the fractional part (on $\MB{C}$)
	\be{align}{eq:fractpart}
	{
		\mbs{z}\coloneqq \mf 12+\mf \iunit{2\pi}\log e^{-\iunit 2\pi \mb{z-\mf 12}}.
	}
	With this definition the function $\mbs{\cdot}$ evaluates on $\MB{R}$ as the usual fractional part. In particular, one verifies quickly
	\[\mbs{z}=\mbs{\Re z}+\iunit \Im z.\]
	On the other hand, setting for every $x\in\MB{R}$,
	\be{align}{eq:sign}
	{
		\on{sgn}\mb{x}\coloneqq\be*{cases}{1&\ft{if }x\ge 1\\0&\ft{if }x< 0},
	}
	we find, using $\Re \mb{1-e^{\iunit 2\pi \on{sgn}\mb{\Im z}z}}\ge 0$ for the splitting of the logarithm,
	\be{align}{eq:logsinfact}
	{
		\log \mb{4\mb{\sin\mb{\pi z}}^2}&=\log \mb{\mb{1-e^{\iunit 2\pi z}}\mb{1-e^{-\iunit 2\pi z}}}\nonumber\\
		&=\log \mb{e^{-\iunit 2\pi \mb{\on{sgn}\mb{\Im z}z-\mf 12}}}+2\log\mb{1-e^{\iunit 2\pi \on{sgn}\mb{\Im z}z}}\\
		=2\pi\iunit&\mb{\mf 12-\mbs{\on{sgn}\mb{\Im z}\Re z+\iunit \mbm{\Im z}}}+2\log\mb{1-e^{\iunit 2\pi \on{sgn}\mb{\Im z}z}}\nonumber.
	}
	Clearly, this last equality signifies the importance to investigate the logarithm of the sine function (log-sine) on $\MB{C}$ as we have exponentially fast convergence to the fractional part in the imaginary part if we take $\Im z\to \pm \infty$.
	For every $s\in\MB{C}$ such that $\Re s>1$ and every $z\in \MB{C}$ such that $\mbm{z}\le 1$ and for every $x\in \MB{C}\setminus\MB{Z}_{\le 0}$, we set
	\be{align}{eq:polylog}
	{
		\operatorname{\Phi}_s\mb{z,x}\coloneqq\sum_{k=1}^\infty\mf{z^k}{\mb{k+x}^s}\ftqq\operatorname{Li}_s\mb{z}\coloneqq \on{\Phi}_s\mb{z,1},		
	}
	and we set for every $\alpha,z\in\MB{C}$ and every $\beta\in\MB{C}\setminus\MB{Z}_{\le 0}$
	\be{align}{eq:confluent}
	{
		{}_1F_1\mb{\alpha,\beta,z}\coloneqq\sum_{k=1}^\infty\mf{\mb{\alpha}_k}{\mb{\beta}_k}\mf {z^k}{k!}.
	}
	
	The following convergence result puts the integral of the fractional part and the log-sine on the same footing:
	\be{theorem}{t:convergence}
	{
		Let $s\in \MB{C}$, $a\in\MB{C}$ such that $\Re a\neq 0$, $b,c\in\MB{C}$ and let $d>0$. Set for every $k\in\MB{N}$
		\be*{align*}
		{
			f_{s,k,1}\mb{a,b,c,d}&\coloneqq \mf{2\pi ad^{1-s}\cos\mb{2\pi\mb{bk+c}}}{s-1}{}_1F_2\mb{\mf 12-\mf s2;\mf 32,\mf 32-\mf s2;-\pi^2a^2d^2k^2},\\
			f_{s,k,2}\mb{a,b,c,d}&\coloneqq \mf{\sin\mb{2\pi\mb{bk+c}}}{sd^sk}{}_1F_2\mb{-\mf s2;\mf 12,\mf 12-\mf s2;-\pi^2a^2d^2k^2},\\
			f_{s,k}\mb{a,b,c,d}&\coloneqq f_{s,k,1}\mb{a,b,c,d}+f_{s,k,2}\mb{a,b,c,d},\\
			g_{s,k}\mb{a,b,c}&\coloneqq k^{s-1}\mb{2\pi}^s\Gamma\mb{-s}\mf 1{2\iunit}\mb{\mb{-\iunit a}^se^{\iunit 2\pi\mb{bk+c}}-\mb{\iunit a}^se^{-\iunit 2\pi\mb{bk+c}}}.
		}
		Set
		\be*{align*}
		{
			d_+\coloneqq\be*{cases}
			{
				-\mf{\Im b}{\Im a},&\substack{\ft{if }\Im a\neq 0\\\ft{and }-\mf{\Im b}{\Im a}>d}\\
				d,&\ft{otherwise}
			}\ftqq d_-\coloneqq\be*{cases}
			{
				1,&\ft{if }\Im a=\Im b=0\\
				\on{sgn}\mb{2d_+\Im a+\Im b},&\ft{otherwise}
			},
		}
		and set
		\be*{align*}
		{
			I_s\mb{a,b,c,d}&\coloneqq\cos\mb{2\pi c}\int_d^\infty\pi\mb{\mf 12-\mbs{ax+b}+\iunit \mb{\Im ax+\Im b}}x^{-s-1}\on{d}x\\
			-&\sin\mb{2\pi c}\int_d^\infty\mb{\mf 12\log\mb{4\mb{\sin\mb{\pi\mb{ax+b}}}^2}-\pi\mbm{\Im ax+\Im b}}x^{-s-1}\on{d}x.
		}
		Then for $\Re s>-1$ it holds
		\be*{align*}
		{
			&I_s\mb{a,b,c,d}=\sum_{k=1}^\infty \cos\mb{2\pi c}\mb{f_{s,k}\mb{\Re a,\Re b,0,d}+g_{s,k}\mb{\Re a,\Re b,0}}\\
			&+\iunit d_-\sin\mb{2\pi c}\mb{f_{s,k}\mb{\Re a,\Re b,0,d}-2f_{s,k}\mb{\Re a,\Re b,0,d_+}-g_{s,k}\mb{\Re a,\Re b,0}}\\
			&+\sin\mb{2\pi c}\mb{\mf{{}_1F_1\mb{-s;1-s;-\iunit 2\pi d_-adk}}{e^{\iunit 2\pi d_-bk}sd^sk}-\mf{{}_1F_1\mb{-s;1-s;-\iunit 2\pi d_-ad_+k}}{e^{\iunit 2\pi d_-bk}sd_+^sk}}\\
			&+\sin\mb{2\pi c}\mb{\mf{{}_1F_1\mb{-s;1-s;\iunit 2\pi d_-ad_+k}}{e^{-\iunit 2\pi d_-bk}sd_+^sk}+g_{s,k}\mb{d_-a,d_-b,\mf 14}+\iunit g_{s,k}\mb{d_-a,d_-b,0}}.
		}
		In particular, if $\Im a=\Im b=0$ it holds
		\[I_s\mb{a,b,c,d}=\sum_{k=1}^\infty f_{s,k}\mb{a,b,c,d}+g_{s,k}\mb{a,b,c}.\]
		The series converges to a holomorphic function in $s$ on $\MB{C}$. The convergence is absolute and uniform subject to the conditions $b,c\in\MB{C}$, $a,d$ bounded away from $0$ and $\infty$ and for every $\epsilon>0$ in $s$ with $\Re s>-\epsilon$ and $\mbm{\Im s}<\epsilon$.
		
		Moreover, the series
		\be{enumerate}{enumroman}
		{
			\item\label{t:convergence:I}$\sum_{k=1}^\infty f_{s,k,1}\mb{\Re a,\Re b,c,d}$,
			\item\label{t:convergence:II}$\sum_{k=1}^\infty f_{s,k,2}\mb{\Re a,\Re b,c,d}$,
			\item\label{t:convergence:III}$\be*{cases}
			{
				\sum_{k=1}^\infty \mf{e^{-\iunit 2\pi d_-bk}}{sk}{}_1F_1\mb{-s;1-s;-\iunit 2\pi d_-axk}&\ft{if }d< x< d_+\\
				\sum_{k=1}^\infty \mf{e^{\iunit 2\pi d_-bk}}{sk}{}_1F_1\mb{-s;1-s;\iunit 2\pi d_-adk}&\ft{if }\Im a= 0}$,
		}
		converge absolutely for every $\Re s<0$ and they satisfy the bounds
		\be{enumerate}{enumroman}
		{
			\setcounter{enumi}{3}
			\item\label{t:convergence:V}$\sum_{k=1}^\infty \mbm{d^{1+s}f_{s,k,1}\mb{\Re a,\Re b,c,d}}\le 2\cosh\mb{\Im c}h_s\mb{\Re a,\Re b,0}$,
			\item\label{t:convergence:VI}$\sum_{k=1}^\infty \mbm{ {d^{s+1}}f_{s,k,2}\mb{\Re a,\Re b,c,d}}\le 2\cosh\mb{\Im c}h_s\mb{\Re a,\Re b,0}$,
			\item\label{t:convergence:VII}$\be*{cases}
			{
				\sum_{k=1}^\infty \mbm{\mf{e^{-\iunit 2\pi d_-bk}}{sk}{}_1F_1\mb{-s;1-s;-\iunit 2\pi d_-axk}}\le h_s\mb{d_-a,d_-b,x}&\ft{if }d< x< d_+\\
				\sum_{k=1}^\infty \mbm{\mf{e^{\iunit 2\pi d_-bk}}{sk}{}_1F_1\mb{-s;1-s;\iunit 2\pi d_-adk}}\le h_s\mb{-d_-a,-d_-b,d}&\ft{if }\Im a= 0
			}$,
		}
		where for some $k_0=k_0\mb{s,a}\in\MB{N}$ set
		\be*{align*}
		{
			h_s\mb{a,b,d}\coloneqq&\mf{\on{Li}_2\mb{e^{2\pi \Im\mb{ad+b}}}}{2\pi\mbm{a}d}+\mf{\on{Li}_2\mb{e^{2\pi k \Im b}}}{2\pi\mbm{a}d}+\mbm{s+1}\mf{\mb{1-e^d}}{d\Re s}\mf{\on{Li}_{1-\Re s}\mb{e^{2\pi k \Im b}}}{\mb{2\pi\mbm{a}}^{-\Re s}}\\
			&+d^{-1}\mb{\mbm{h_{s,k_0}\mb{a,b}}+\mbm{h_{s,k_0}\mb{a,b}-h_{s,1}\mb{a,b}}},\\
			h_{s,k}\mb{a,b}\coloneqq& 2\mb{s+1}\mb{\mf{\Phi_2\mb{e^{2\pi\Im b},k}}{2\pi\mbm{a}\mb{1+\Re s}}-\mf{\Phi_{1-\Re s}\mb{e^{2\pi\Im b},k}}{\mb{2\pi \mbm{a}}^{-\Re s}\mb{1+\Re s}}},
		}
		and for every $s\in\MB{C}$ it holds
		\be{enumerate}{enumroman}
		{
			\setcounter{enumi}{6}
			\item\label{t:convergence:VIII}
			\be*{align*}
			{
				\sum_{k=1}^\infty&\mbm{\mf{e^{-\iunit 2\pi d_-bk}}{sk}\mb{\mf{{}_1F_1\mb{-s;1-s;-\iunit 2\pi d_-adk}}{d^s}-\mf{{}_1F_1\mb{-s;1-s;-\iunit 2\pi d_-ad_+k}}{d_+^s}}}\le\\
				& \mf {\zeta\mb{2}}{2\pi\mbm{a}}\mb{\sup_{k\in\MB{N}}\mbm{\mf 1{d_+^{ s+1}}-\mf {e^{-\iunit 2\pi d_-a\mb{d-d_+}k}}{d^{s+1}}}+\mf {\mbm{s+1}}{\mbm{1+\Re s}}\mbm{d^{-s-1}-d_+^{-s-1}}},
			}
			\item\label{t:convergence:IX}for every $n\in\MB{N}$ such that $\Re s+n>0$
			\be*{align*}
			{
				\sum_{k=1}^\infty&\mbm{\mf{{}_1F_1\mb{-s;1-s;\iunit 2\pi d_-ad_+k}}{e^{-\iunit 2\pi d_-bk}sd_+^sk}+g_{s,k}\mb{d_-a,d_-b,\mf 14}+\iunit g_{s,k}\mb{d_-a,d_-b,0}}\le\\
				&\sum_{m=0}^{n-1}\mf{\mbm{\mb{s+1}_m}\on{Li}_{m+1}\mb{e^{-2\pi d_-\Im\mb{ad_++b}}}}{\mb{2\pi\mbm{a}}^{m+1}d^{\Re s+m+1}}+\mf{\mbm{\mb{s+1}_n}\on{Li}_{n+1}\mb{e^{-2\pi d_-\Im\mb{ad_++b}}}}{\mb{2\pi \mbm{a}}^nd^{\Re s+n}\mb{\Re s+n}}.
			}
		}
	}
	Concerning the general validity of the above result, consider a function $F$ with poles and zeros of order in $\MB{C}$ such that for almost every $x\in\MB{R}$
	\[\log F\mb{x}=\sum_{k=1}^\infty a_ke^{\iunit 2\pi xb_k},\]
	for some $a_k,b_k\in\MB{C}$ (e.g. functions of the form $\prod_{k=1}^\infty\mb{\sin\mb{b_kx}e^{-\iunit\pi b_kx }2\iunit}^{a_k}$ which are closely connected to integrals of Lambert series). Then the integral
	\[\int_d^\infty \log F\mb{x}x^{-s-1}\on{d}x\]
	would follow a similar evaluation, ignoring convergence issues here, showing that the expansion in ${}_1F_2$ is intimately connected to the represetation of $\log F$ as an exponential series. Conversely, since $F$ has a specific zero and pole structure, the alternative evaluation analogous to the \lhs from \tm{maintheorem:B} would pick up the poles and zeroes of $F$, resulting in a sum in ${}_2F_1$ evaluated at these points. Thus this gives the general identity in \tm{maintheorem:B} somehow the nature of a spectral sum. Moreover, this identifies series in $k^{-s}$  as rather special functions, since as we will see in the proofs for \tm{convergence} and \tm{maintheorem:B} in \se{proof}, on the \lhs $k^{-s}$ emerges from the zeroes of the sine-function, involving the identity in \eq*{hyperidentity}, and on the \rhs terms with factors $k^{s-1}$ are a consequence of the Fourier modes, involving the identity in \eq*{gammaintegral}.\\
	In what follows we discuss the integrand in $I_s$, assuming the integral kernel to be  $x^{-s-1}$ which can be seen as the canonical extension for the radial measure to arbitrary complex dimension. In fact, the choice of parameters is natural in the following sense. Set
	\[h\mb{x}\coloneqq {1-e^{\iunit 2\pi x}}=\iunit 2\sin\mb{\pi x}e^{\iunit\pi x}.\]
	Then we have the identities
	\be*{align*}
	{
		h_-\mb{x}&\coloneqq\log\mf{h\mb{x}}{h\mb{-x}}=-\iunit 2\pi\mb{\mf12-\mbs{x}}\ftqq\\  h_+\mb{x}&\coloneqq\log\mb{h\mb{x}h\mb{-x}}=\log \mb{4\mb{\sin\mb{\pi z}}^2}.
	}
	and the composition of $h$ with an affine complex map yields precisely the parameters $a,b\in\MB{C}$ and the parameter $c$ is just an interpolation parameter between the two functions $h_+$ and $h_-$ while the parameter $d$ serves as a cut-off parameter which plays also a role in the analysis as can be seen by the appearance of $d_+$ in place of $d$. In what follows, we consider different (not necessarily independent) point of views concerning the nature of the integrand $h_\pm$.

	\subsection{Analytic interpretation}
	
	We may see the functions inside the $\log$ of $h_+$ and $h_-$ as a geometric section for $\log \mb{1-e^{\iunit 2\pi x}}$ for the modulus $2$. What we mean by geoemtric section is the product analogue of an arithmetic section: if we set for every $n\in\MB{N}$ and every function $f\colon \MB{C}\rightarrow\MB{C}$
	\be*{align*}
	{
		\omega_n\coloneqq \mb{-1}^{\mf 2n}\ftqq\on{GS}_{n,\ell}\mb{f}\coloneqq \prod_{k=0}^{n-1}\mb{f\circ \omega_n^k\on{id}}^{\omega_n^{kl}},
	}
	then $h_{\pm}=\log \on{GS}_{2,\mf{1\mp 1}2}\mb{1-e^{\iunit 2\pi \cdot}}$. This point of view gives then an immediate generalization of the integrand to geometric sections for higher modulus.	Moreover, the following corollary signifies that the composition by a complex affine map is not an arbitrary generalization:
	
	\be{corollary}{c:limit}
	{
		Let $s\in\MB{C}$, $d>0$ and let $a,b\in\MB{C}$, $a\neq 0$ such that $d_+=d$ and $d_-=1$, recalling \tm{convergence}. Denoting the holomorphic continuation of $I_s\mb{a,b,c,d}$ in $s$ by $\widetilde{I}_s\mb{a,b,c,d}$, then it holds exponentially fast
		\be*{align*}
		{
			\lim_{d\Im a+\Im b\to\infty}\widetilde{I}_s\mb{a,b,\mf 14,d}=-\iunit \widetilde{I}_s\mb{a,b,0,d}.
		}
	}
	\be*{proof}
	{
		By assumption $d_+=d$ and $d_-=1$. Then, considering the bound in \tm[IX]{convergence} which decreases exponentially fast to $0$, it follows
		\be*{align*}
		{
			\lim_{d\Im a+\Im b\to\infty}\iunit\widetilde{I}_s\mb{a,b,\mf 14,d}=\sum_{k=1}^\infty f_{s,k}\mb{\Re a,\Re b,0,d}+g_k\mb{\Re a,\Re b,0}=\widetilde{I}_s\mb{a,b,0,d}.
		}
	}
	
	\subsection{Geometric interpretation}
	
	The set $S\coloneqq \mbs{z\in\MB{C}\SR \mbm{z}=1}$ is uniquely linearly parameterized by $e^{\iunit 2\pi x}$ for every $x\in\mblo{-\mf 12,\mf12}$. The direction in $\MB{C}$ to the real axis at $1$ is
	\[1-e^{\iunit 2\pi x}=h\mb{x}.\]
	This is also the shortest direction between these two points if one is allowed to pass through the unit disc. The quadratic length of this vector is
	\[\mbm{h\mb{x}}^2=\mbm{1-e^{\iunit 2\pi x}}^2=\mb{1-e^{\iunit 2\pi x}}\mb{1-e^{-\iunit 2\pi x}}=4\mb{\sin\mb{\pi x}}^2=e^{h_+\mb{x}},\]
	which is of course also the squared length of the corresponding line segment in the Eucledian plane $\MB{R}^2$.
	On the other hand, the corresponding shortest direction on the unit circle from $1$ to $e^{\iunit 2\pi x}$ is given by
	\[2\pi x=2\pi\mb{\mbs{x+\mf 12}-\mf 12}\ftqq[if]x\not\in \mbs{-\mf 12,\mf 12}.\]
	Now the quotient of the length through the disc with the length on the circle is given for $-\mf 12<x<\mf 12$ by
	\[\mf{2\mbm{\sin\mb{\pi x}}}{\mbm{2\pi x}}=\mf{\sin\mb{\pi x}}{\pi x}.\]
	However in the integrand for $c=\mf 14$ we take the logarithm of $e^{h_+\mb{x}}$ which gives a rather unintuitve interpretation. For once for small values in $x$, using Taylor expansion, we find 
	\[\mf{\log\mbm{\sin\mb{\pi x}}}{\mbm{\pi x}}\sim \mf{\log\mbm{\pi x}}{\mbm{\pi x}}.\]
	Since $\log\mbm{\sin\mb{\pi x}}$ is unbounded we cannot interpret $h_+\mb{x}$ again as the shortest distance on a corresponding circle but rather as a weight for a shortest Eulidean path through the unit disc. So the case $c=0$ and $c=\mf 14$ seem to be of significant different geometric origin. Regarding the values $\zeta\mb{2n}$ and $\zeta\mb{2n+1}$ and their relation to the case $c=0$ and $c=\mf 14$ in view of \lem{recurrence:odd}-\lem*{recurrence:even}, respectively, it seems therefore rather unlikely that $\zeta\mb{2n+1}\pi^{-2n-1}$ is a rational number. Moreover, since one can pass from $\zeta\mb{s}$ to $\zeta\mb{s+1}$ holomorphically there seems to be captured a repeated interpolation between the two geometric origins within the zeta function. Assuming interpolation with constant spacing, i.e. attaching an open strip $\mb{n-\mf 12,n+\mf12}+\iunit\MB{R}$ to each integer $n\in\MB{Z}$, this suggests that a result of the form 
	\[\forall s\in{I}\colon\zeta\mb{s}\pi^{-s}\in\MB{Q} \Rightarrow \zeta\mb{s+1}\pi^{-s-1}\not\in\MB{Q},\]
	should hold for some set $I\subseteq\MB{R}$ containing $\MB{Z}$ as large as maybe $\MB{R}\setminus\mb{\MB{Z}+\mf 12}$. Note that this is true for the integers $I=\MB{Z}_{\le 1}$.
	
	\subsection{Dynamical interpretation}
	Among the arguably most simple dynamical systems is the continuous rotation of a single particle on a circle at constant speed. Closely related is the billiard for one particle (we assume specular reflection at the boundary) with a circle as boundary. While this billiard is a continuous dynamical system the evolution of the intersection of the particle with the boundary forms a discrete dynamical system, namely a discrete rotation on the unit circle. Now the evolution of the latter is given by 
	\[S\times S\ni\mb{e^{\iunit 2\pi x},e^{\iunit 2\pi y}}\mapsto \mb{e^{\iunit 2\pi x+y},e^{\iunit 2\pi y}}\in S\times S.\]
	More precisely the set $S\times S$ describes the tangent space (in physical terms the phase space) for the billiard at the boundary in some coordinate system. We notice the rotational symmetry in the first component since the reflection at the boundary is independent of the place on the boundary but depends only on the condition for the initial direction for the particle. Hence the reduction by rotational symmetry in the first component yields a dynamical system with exact localization in positional space with the trivial evolution
	\[\mbs{1}\times S\ni\mb{1,e^{\iunit 2\pi y}}\mapsto \mb{1,e^{\iunit 2\pi y}}\in \mbs{1}\times S,\]
	where we took $1\in S$ as the representative for $S$. And while this system is trivial it makes sense as an accelerated system on periodic orbits (for the billiard or at the boundary), namely whenever $y\in\MB{Q}\cap \mblo{-\mf 12,\mf 12}$ since on these orbits the particle will return to $1$. This would be in fact a valid game strategy for a pinball with one pad, where the ball is shot in the direction of a periodic orbit since then it is sure that the ball will never get stuck. The pinball we derived from the billiard we started here with would be the most boring realization of one. Moreover, for a somewhat ideal pinball the only configuration which matters is the positioning of the pads which is in turn for our pinball related to
	\[1-e^{\iunit 2\pi x}=h\mb{x},\]
	which in turn can be identified with an (equivalence class of an) orbit for the underlying billiard or for the underlying discrete dynamical system on the boundary. To this end, it should be noted that this pinball distinguishes the periodic orbits as a set up to mirror symmetry while the system on the boundary makes no distinction at all. In particular, a primitive periodic orbit at the boundary of length $n\in\MB{N}$ has a degeneracy of $\varphi\mb{n}$ (where $\varphi$ denotes the totient function). Notably, by comparing the primitive periodic orbits with these for the corresponding billiard this leads to an immediate dynamical proof that $\varphi\mb{n}$ is an odd natural number if and only if $n\in\mbs{1,2}$. Moreover, the collection of primitive periodic orbits on the boundary up to length $n$ is parameterized by the Farey sequence of order $n$ and a precise description how this sequence distributes on the unit circle as $n\to\infty$ is ongoing research to which belongs an open question equivalent to the Riemann hypothesis \cite{Franel_1924}. In the context of the Riemann hypothesis, a connetion to escape rates for open circular open billiards has been made \cite{Bunimovich_2005} and their setting suggests an immediate connection to the concept of hohlraum radiation in physics where we would interpret the escape of orbits though a hole as leaking radiation (physically) measured by an outer observer. 
	
	\section{Proof of \tm{convergence} and \tm{maintheorem}-\tm*{maintheorem:B}}\label{s:proof}
	For every $a,b,c\in\MB{C}$, every $d>0$ and every $s\in\MB{C}$ such that $s\not\in\MB{Z}_{\le 0}$ we set
	
	\be{align}{eq:fk}
	{
		f_{s,k}&\mb{a,b,c,d}\coloneqq\mf{\sin\mb{2\pi\mb{bk+c}}}{d^ssk}{}_1F_2\mb{-\mf s2;\mf 12,1-\mf s2;-\pi^2a^2d^2k^2}\nonumber\\
		&-\mf{2\pi a\cos\mb{2\pi\mb{bk+c}}d^{1-s}}{1-s}{}_1F_2\mb{\mf 12\mb{1-s};\mf 32,1+\mf 12\mb{1-s};-\pi^2a^2d^2k^2}.
	}
	
	We recall ${}_1F_1$ defined in \eq*{confluent}.
	
	\be{lemma}{l:mellintransform:I}
	{
		Let $k\in\MB{N}$, $s\in\MB{C}$ such that $\Re s<0$, let $a,b,c\in\MB{C}$ and let $d>0$. Then it holds
		\be{align}{eq:sinintegral}
		{
			f_{s,k}\mb{a,b,c,d}=-\mf 1k\int_0^d\sin\mb{2\pi \mb{akx+bk+c}}x^{-s-1}\on{d}x.
		}
		In particular, for every $k\in\MB{N}$ the above integral admits a meromorphic continuation in $s\in\MB{C}$ given by $f_{s,k}$ with (possible) poles at $s\in\MB{Z}_{\ge 0}$ and a holomorphic continuation in $d\in\MB{C}\setminus\MB{R}_{\le 0}$ both holomorphic in $a,b,c\in\MB{C}$.
		Moreover, it holds
		\be{align}{eq:expidentity}
		{
			f_{s,k}\mb{a,b,c+\mf 14,d}+\iunit f_{s,k}\mb{a,b,c,d}=\mf{e^{\iunit 2\pi 	\mb{bk+c}}}{sd^sk}{}_1F_1\mb{-s;1-s;\iunit 2\pi adk}.
		}
	}
	\be*{proof}
	{
		We recall the general Euler transform \cite[Equation (4.1.2)]{Slater_1969} and \cite[Sect. 3.6, Equation (10)]{Luke_1969} which in our case simplifies for every $\Re s<1$, for every $\beta\in\MB{C}\setminus\MB{Z}_{\le 0}$ and every $z\in\MB{C}$ to
		\be{align}{eq:EulerTransform}
		{
			{}_1F_2\mb{\mf 12-\mf s2;\beta, \mf32-\mf s2;z}=\mf{1-s}2\int_0^1{}_0F_1\mb{\beta;xz}x^{-\mf 12-\mf s2}\on{d}x,
		}
		where we set (analogously to the definition of ${}_1F_2$ in \eq*{hypergeo12})
		\be*{align*}
		{
			{}_0F_1\mb{\beta;z}\coloneqq \sum_{n=0}^\infty \mf 1{\mb{\beta}_n}\mf{z^n}{n!}.
		}
		Comparing Taylor series, we find the identities
		\[\mf{\sin x}x={}_0F_1\mb{\frac 32;-\mf{x^2}4}\ftqq\cos x={}_0F_1\mb{\frac 12;-\mf{x^2}4}.\]
		We let $\Re s<1$ and $\alpha>0$ and compute, using in addition \eq*{EulerTransform} in which we take $\beta=\mf 32$,
		\be{align}{eq:sinidentity}
		{
			\int_0^d\sin\mb{\alpha x}x^{-s-1}\on{d}x&=\alpha d^{1-s}\int_0^1{}_0F_1\mb{\frac 32;-\mf{\alpha^2d^2}4 x^2}x^{-s}\on{d}x\nonumber\\
			&= \alpha d^{1-s} 2\int_0^1{}_0F_1\mb{\frac 32;-\mf{\alpha^2d^2}4 x}x^{-\mf {s+1}2}\on{d}x\nonumber\\
			&=\mf{\alpha d^{1-s}}{1-s}{}_1F_2\mb{\mf 12\mb{1-s};\frac 32,1+\mf 12\mb{1-s};-\mf{\alpha^2d^2}4}.		
		}
		On the other hand, we let $\Re s<0$ and compute, using in addition \eq*{EulerTransform} in which we take $\beta=\mf 12$ and $s\mapsto s+1$,
		\be{align}{eq:cosidentity}
		{
			\int_0^d\cos\mb{\alpha x}x^{-s-1}\on{d}x&=d^{-s}\int_0^1{}_0F_1\mb{\frac 12;-\mf{\alpha^2d^2}4 x^2}x^{-s-1}\on{d}x\nonumber\\
			&=\mf 12d^{-s}\int_0^1{}_0F_1\mb{\frac 12;-\mf{\alpha^2d^2}4 x}x^{-\mf s2}\on{d}x\nonumber\\
			&=-\mf {d^{-s}}s{}_1F_2\mb{-\mf s2;\frac 12,1-\mf s2;-\mf{\alpha^2d^2}4}.		
		}
		We note for every $\beta\in\MB{C}$
		\be{align}{eq:exponentialexpansion}
		{
			e^{\iunit 2\pi\mb{\beta+bk+c}}&=\cos\mb{2\pi\mb{bk+c}}\cos\mb{2\pi\beta}-\sin\mb{2\pi\mb{bk+c}}\sin\mb{2\pi\beta}\nonumber\\
			+&\iunit \cos\mb{2\pi\mb{bk+c}}\sin\mb{2\pi\beta}+\iunit\sin\mb{2\pi\mb{bk+c}}\cos\mb{2\pi\beta}.
		}
		Taking now $\alpha=2\pi ak$ in \eq*{sinidentity}-\eq*{cosidentity}, and combining it with the identity in \eq*{exponentialexpansion} in which we take $\beta=2\pi akx$, we find
		\be*{align*}
		{
			f_{s,k}&\mb{a,b,c,d}=\mf{\sin\mb{2\pi\mb{bk+c}}}{d^ss}{}_1F_2\mb{-\mf s2;\mf 12,1-\mf s2;-\pi^2a^2d^2k^2}\\
			&-\mf{2\pi akd^{1-s}\cos\mb{2\pi\mb{bk+c}}}{1-s}{}_1F_2\mb{\mf 12\mb{1-s};\mf 32,1+\mf 12\mb{1-s};-\pi^2a^2d^2k^2},
		}
		and we conclude, recalling the possible poles for ${}_1F_2$ (see \eq*{hypergeo12}). The identity in \eq*{expidentity}, using again an analogous transformation as in \eq*{EulerTransform} with $\beta$ removed, noting
		\[{}_0F_0\mb{z}=e^z.\] 
	}
	
	\be{lemma}{l:mellintransform:I:recurrence}
	{
		Let $n,k\in\MB{N}$, $s\in\MB{C}\setminus \MB{N}$, $a,d\in\MB{C}\setminus\mbs{0}$. Then it holds
		\be*{align*}
		{
			{}_1F_1\mb{-s;1-s;\iunit 2\pi adk}=&-e^{\iunit 2\pi kad}\sum_{m=1}^n\mf{\mb{s}_m}{\mb{\iunit 2\pi kad}^m}\\
			&+\mf{\mb{s}_n}{\mb{\iunit 2\pi kad}^n}{}_1F_1\mb{-s-n;1-s-n;\iunit 2\pi adk},
		}
		and if $\Re s<0$ it holds
		\be{align}{eq:mellintransform:I:bound:I}
		{
			&\mf{\mbm{{}_1F_1\mb{-s-1;1-s-1;\iunit 2\pi adk}}}{\mbm{s+1}}\le \mf 1{\mbm{s+1}}\nonumber\\
			&+\mb{2\pi \mbm{a}k}^{\Re s+1}\mb{\mf{1-e^d}{\Re s}+\sup_{x\in\mb{d, 2\pi \mbm{a}dk}}\mbm{e^{\iunit \mf a{\mbm{a}}x}-1}\mbm{\mf{\mb{2\pi \mbm{a}k}^{-\Re s-1}-1}{1+\Re s}}},
		}
		and if $s\in\MB{C}$ and $d_1,d_2>0$ it holds
		\be{align}{eq:mellintransform:I:bound:II}
		{
			\mf{2\pi\mbm{a}k}{\mbm{s}}&\mbm{\mf{{}_1F_1\mb{-s;1-s;\iunit 2\pi ad_1k}}{d_1^s}-\mf{{}_1F_1\mb{-s;1-s;\iunit 2\pi ad_2k}}{d_2^s}}\le\nonumber\\
			&{\mbm{\mf{e^{\iunit 2\pi ad_1k}}{d_1^{s+1}}-\mf{e^{\iunit 2\pi ad_2k}}{d_2^{s+1}}}+\mf{\mbm{s+1}\mbm{d_1^{-\Re s-1}-d_2^{-\Re s-1}}}{\mbm{1+\Re s}}\sup_{x\in\mbs{d_1,d_2}}e^{-2\pi kx\Im a}}.
		}
	}
	\be*{proof}
	{
		We compute, using the identity in \eq*{expidentity} from \lem{mellintransform:I} and integration by parts, and assuming $\Re s<-1$,
		\be*{align*}
		{
			&\mf 1{sd^s}{}_1F_1\mb{-s;1-s;\iunit 2\pi adk}=kf_{s,k}\mb{a,0,\mf 14,d}+\iunit kf_{s,k}\mb{a,0,0,d}\\
			&=-\int_0^d e^{\iunit 2\pi kax}x^{-s-1}\on{d}x
			=-\mf{d^{-s-1}e^{\iunit 2\pi kad}}{\iunit 2\pi ka}-\mf{s+1}{\iunit 2\pi ka}\int_0^d e^{\iunit 2\pi kax}x^{-s-2}\on{d}x\\
			&=-\mf{e^{\iunit 2\pi kad}}{\iunit 2\pi kad^{s+1}}+\mf 1{\iunit 2\pi kad^{s+1}}{}_1F_1\mb{-\mb{s+1};1-\mb{s+1};\iunit 2\pi adk}.
		}
		Hence by meromorphic continuation this recurrence relation holds for all $s\in\MB{C}$ away from the poles for ${}_1F_1$. Then inducing on this recurrence relation, this yields the claim on the recurrence relation. To see the bound in \eq*{mellintransform:I:bound:I}, we estimate for $\Re s<0$
		\be*{align*}
		{
			&\mbm{\int_0^d \mb{e^{\iunit 2\pi kax}-1}x^{-s-2}\on{d}x}=\mbm{\mb{2\pi \mbm{a}k}^{s+1}\mb{\int_0^d+\int_d^{2\pi \mbm{a}dk}} \mf{e^{\iunit \mf a{\mbm{a}}x}-1}{x^{s+2}}\on{d}x}\\
			&\le\mf{\mb{2\pi \mbm{a}k}^{\Re s+1}}{d^{\Re s +1}}\mb{\mf{1-e^d}{\Re s}+\sup_{x\in\mb{d, 2\pi \mbm{a}dk}}\mbm{e^{\iunit \mf a{\mbm{a}}x}-1}\mbm{\mf{\mb{2\pi \mbm{a}k}^{-\Re s-1}-1}{1+\Re s}}},
		}
		and in particular the integral on the \lhs is holomorphic in $s$ such that $\Re s<0$, and we calculate, assuming $\Re s<-1$,
		\[\int_0^d x^{-s-2}\on{d}x=-\mf{d^{-s-1}}{s+1},\]
		which has a meromorphic continuation in $s\in\MB{C}$.
		To see the bound in \eq*{mellintransform:I:bound:II}, we estimate, integrating by parts,
		\be*{align*}
		{
			\mf{2\pi\mbm{a}k}{\mbm{s}}&\mbm{d_1^{-s}{}_1F_1\mb{-s;1-s;\iunit 2\pi ad_1k}-d_2^{-s}{}_1F_1\mb{-s;1-s;\iunit 2\pi ad_2k}}\\
			&\le {\mbm{\mf{e^{\iunit 2\pi ad_1k}}{d_1^{s+1}}-\mf{e^{\iunit 2\pi ad_2k}}{d_2^{s+1}}}+\mbm{s+1}\mbm{\int_{d_1}^{d_2}e^{\iunit 2\pi kax}x^{-s-2}\on{d}x}}\\
			&\le {\mbm{\mf{e^{\iunit 2\pi ad_1k}}{d_1^{s+1}}-\mf{e^{\iunit 2\pi ad_2k}}{d_2^{s+1}}}+\mf{\mbm{s+1}\mbm{d_1^{-\Re s-1}-d_2^{-\Re s-1}}}{\mbm{1+\Re s}}\sup_{x\in\mbs{d_1,d_2}}e^{-2\pi kx\Im a}}.
		}
	}
	
	Moreover, for every $k\in\MB{N}$, $s\in\MB{C}\setminus\MB{Z}_{\ge 0}$, $a,b,c\in\MB{C}$, $a\neq 0$, $d>0$ we set
	\be*{align}
	{
		g_{s,k}\mb{a,b,c}\coloneqq  \mf {k^{s-1}}{2\iunit}\mb{2\pi}^s\Gamma\mb{-s}\mb{e^{\iunit 2\pi\mb{bk+c}}\mb{-\iunit a}^s-e^{-\iunit 2\pi\mb{bk+c}}\mb{\iunit a}^s},
	}
	and we set
	\be{align}{eq:Fk}
	{
		F_{s,k}\mb{a,b,c,d}&\coloneqq f_{s,k}\mb{a,b,c,d}+g_{s,k}\mb{a,b,c}.
	}
	
	\be{lemma}{l:mellintransform:II}
	{
		Let $s\in\MB{C}$ such that $\Re s>-1, k\in\MB{N}$, $a\in \MB{R}\setminus\mbs{0}$, $b,c\in\MB{C}$ and let $d>0$. Then it holds
		\[F_{s,k}\mb{a,b,c,d}=k^{-1}\int_d^\infty\sin\mb{2\pi\mb{axk+bk+c}}x^{-s-1}\on{d}x.\]
		In particular, for every $b,c\in\MB{C}$ and every $a\in\MB{R}\setminus\mbs{0}$ and every $d>0$ the function $F_{s,k}\mb{a,b,c,d}$ is holomorphic in $s\in\MB{C}$ and $C^\infty$ in $a$  and if $s\not\in\MB{Z}_{\le 0}$ it is holomorphic in $a\in\MB{C}\setminus\iunit \MB{R}$ and in $d\in\MB{C}\setminus\MB{R}_{\le 0}$, and, moreover,
		\[F_{s,k}\mb{a,b,c+\mf 14,d}+\iunit F_{s,k}\mb{a,b,c,d}\]
		admits a holomorphic continuation in $a\in\MB{C}\setminus\iunit \MB{R}_{\le 0}$ if $s\in\MB{C}\setminus\MB{Z}_{\le 0}$.
	}
	\be*{proof}
	{
		We write, using \lem{mellintransform:I}, assuming first $-1<\Re s <0$,
		\be*{align*}	
		{
			\int_d^\infty\sin\mb{2\pi\mb{axk+bk+c}}x^{-s-1}\on{d}x&=\int_0^\infty\sin\mb{2\pi\mb{axk+bk+c}}x^{-s-1}\on{d}x\\
			+&f_{s,k}\mb{a,b,c,d}.
		}
		We recall the identity for $-1<\Re s <0$ (a proof in the case $a=-1$ is given in \cite[p.51, Satz 2.21]{Remmert_2007} or \cite[p.164]{Baker_2012} and for general $a\in \MB{R}\setminus\mbs{0}$ use substitution and symmetry of the cosinus):
		\be{align}{eq:Gamma}
		{
			\int_0^\infty e^{\iunit ax}x^{-s-1}\on{d}x=\mb{-\iunit a}^s\Gamma\mb{-s}.
		}
		Note that the \rhs in \eq*{Gamma} is holomorphic in $-1<\Re s<0$ and also is $f_{s,k}$, recalling \lem{mellintransform:I}, and extends holomorphically in $a\in \MB{C}\setminus\iunit \MB{R}_{\le 0}$, recalling the logarithm and \eq*{complexpow}. The final claim, including the statements on $b,c,d$, follows immediately, using the identity
		\[2\iunit\sin\mb{2\pi\mb{ax+b+c}}=e^{\iunit 2\pi\mb{ax+b+c}}-e^{-\iunit 2\pi\mb{ax+b+c}},\]
		together with \lem{mellintransform:I}. To see that the claimed integral identity extends in $\Re s>-1$, we integrate by parts, assuming $\Re s>0$ and $\Im a\ge 0$,
		\be*{align*}
		{
			\int_d^\infty e^{\iunit 2\pi axk}x^{-s-1}\on{d}x=-\mf{e^{\iunit 2\pi adk}}{\iunit 2\pi akd^{s+1}}+\mf{s+1}{\iunit 2\pi ak}\int_d^\infty e^{\iunit 2\pi axk}x^{-s-2}\on{d}x,
		}
		and we see that the \rhs extends indeed holomorphically in $\Re s>-1$ if $\Im a\ge 0$. Hence this extension must coincide with the \rhs in \eq*{Gamma} and we conclude.
	}
	
	\be{lemma}{l:mellintransform:II:recurrence}
	{
		For every $n\in\MB{N}$, for every $s\in\MB{C}$, for every $a,d\in\MB{C}\setminus\mbs{0}$ and every $b,c\in\MB{C}$ it holds
		\be*{align*}
		{
			F_{s,k}&\mb{a,b,c+\mf 14,d}+\iunit F_{s,k}\mb{a,b,c,d}=-\sum_{m=0}^{n-1}\mb{s+1}_m\mf{e^{\iunit 2\pi\mb{adk+bk+c}}}{\mb{\iunit 2\pi adk}^{m+1} kd^s}\\
			&+\mf {\mb{s+1}_n}{\mb{\iunit 2\pi ak}^n}\mb{F_{s+n,k}\mb{a,b,c+\mf 14,d}+\iunit F_{s+n,k}\mb{a,b,c,d}}.
		}
		If $\Re s>0$, $\Im a\ge 0$ and $d>0$, it holds
		\be*{align*}
		{
			\mbm{F_{s,k}\mb{a,b,c+\mf 14,d}+\iunit F_{s,k}\mb{a,b,c,d}}\le e^{-2\pi\Im\mb{adk+bk+c}}\mf {d^{-\Re s}}{k\Re s}.
		}		
	}
	\be*{proof}
	{
		We calculate, using \lem{mellintransform:II}, assuming first $-1<\Re s<0$, $\Im a\ge 0$, $d>0$, then using integration by parts, noting that the second equality is valid for $\Re s>-1$ if $\Im a\ge 0$,
		\be*{align}
		{
			kF_{s,k}&\mb{s,a,b,c+\mf 14,d}+k\iunit F_{s,k}\mb{a,b,c,d}=\int_d^\infty e^{\iunit 2\pi\mb{axk+bk+c}}x^{-s-1}\on{d}x\nonumber\\
			&=-\mf{e^{\iunit 2\pi\mb{adk+bk+c}}}{\iunit 2\pi a d^{1+s}k}+\mf {\mb{s+1}}{\iunit 2\pi ak}\int_d^\infty e^{\iunit 2\pi\mb{axk+bk+c}}x^{-s-2}\on{d}x\nonumber\\
			\label{eq:partialintegration}&=-\mf{e^{\iunit 2\pi\mb{adk+bk+c}}}{\iunit 2\pi a d^{1+s}k}+\mf {\mb{s+1}}{\iunit 2\pi a}\mb{F_{s+1,k}\mb{a,b,c+\mf 14,d}+\iunit F_{s+1,k}\mb{a,b,c,d}}.
		}
		Now by (meromorphic) continuation, using \lem{mellintransform:II}, the equality holds on the prescribed domains. Inducing on the equality in \eq*{partialintegration} $n-1$ times, where $n\in\MB{N}$, we conclude. To see the bound, we bound trivially the \rhs in the first equality above
		\[\mbm{\int_d^\infty e^{\iunit 2\pi\mb{axk+bk+c}}x^{-s-1}\on{d}x}\le e^{-2\pi\Im\mb{adk+bk+c}}\mf {d^{-\Re s}}{\Re s}.\]
	}
	We have now finished the preparations for our proof for \tm{convergence}:
	\be*[Proof of \tm{convergence}]{proof}
	{
		We set
		\be*{align*}
		{
			d_+&\coloneqq\be*{cases}
			{
				-\mf{\Im b}{\Im a},&\ft{if }\Im a\neq 0\ft{ and }-\mf{\Im b}{\Im a}>d\\
				d,&\ft{otherwise}
			},\ft{ and}\\
			d_-&\coloneqq\be*{cases}
			{
				1,&\ft{if }\Im a=0=\Im b\\
				\on{sgn}\mb{\Im b}&\ft{if }\Im a=0\neq\Im b\\
				\lim_{\epsilon\to 0} \on{sgn}\mb{\Im a\mb{d_++\epsilon}+\Im b},&\ft{if }\Im a\neq 0
			}.
		}
		Then we write, recalling the factorization in \eq*{logsinfact} and $\on{sgn}$ in \eq*{sign},
		\be*{align*}
		{
			&\mf 12\int_d^\infty\mb{\log\mb{4\mb{\sin\mb{\pi\mb{ax+b}}}^2}-\pi\mbm{x\Im a+\Im b}}x^{-s-1}\on{d}x=\\
			&\int_d^\infty\mb{\iunit\pi \on{sgn}\mb{x\Im a+\Im b}\mb{\mf 12-\mbs{x\Re a+\Re b}}+\log\mb{1-e^{\iunit 2\pi \on{sgn}\mb{x\Im a+\Im b}\mb{ax+b}}}}\mf{\on{d}x}{x^{s+1}}\\
			&=\iunit \pi d_-\mb{\int_{d_+}^\infty-\int_0^{d_+}+\int_0^d}\mb{\mf 12-\mbs{x\Re a+\Re b}}x^{-s-1}\on{d}x\\
			&+\mb{\int_0^{d_+}-\int_0^d}\log\mb{1-e^{-\iunit 2\pi d_-\mb{ax+b}}}\mf{\on{d}x}{x^{s+1}}+\int_{d_+}^\infty\log\mb{1-e^{\iunit 2\pi d_- \mb{ax+b}}}\mf{\on{d}x}{x^{s+1}}.
		}
		By assumption $d_-\Im\mb{ax+b}\le 0$ for every $x\in\mb{0, d_+}$. Hence, using the Taylor series for $\log\mb{1-\cdot}$, recalling $f_{s,k}$ defined in \eq*{fk}, we compute, using \eq*{sinintegral}-\eq*{expidentity} from \lem{mellintransform:I}, assuming $\Re s<0$,
		\be*{align*}
		{
			&\mb{\int_0^{d_+}-\int_0^d}\log\mb{1-e^{-\iunit 2\pi d_-\mb{ax+b}}}x^{-s-1}\on{d}x\\
			&=\sum_{k=1}^\infty f_{s,k}\mb{d_-a,d_-b,\mf 14,d_+}-f_{s,k}\mb{d_-a,d_-b,\mf 14,d}\\
			&\quad\quad\quad\quad\quad\quad\quad\quad+\iunit \mb{f_{s,k}\mb{d_-a,d_-b,0, d}-f_{s,k}\mb{d_-a,d_-b,0, d_+}}\\
			&=\sum_{k=1}^\infty \mf{e^{-\iunit 2\pi d_-bk}}{sk}\mb{d_+^{-s}{}_1F_1\mb{-s;1-s;-\iunit 2\pi d_-ad_+k}-d^{-s}{}_1F_1\mb{-s;1-s;-\iunit 2\pi d_-adk}}.
		}
		On the other hand, using \lem{mellintransform:II} and, recalling the definition of $F_{s,k}$ in \eq*{Fk}, the identities in \eq*{sinintegral}-\eq*{expidentity} from \lem{mellintransform:I}, we find, assuming $\Re s>-1$,
		\be*{align*}
		{
			&\int_{d_+}^\infty\log\mb{1-e^{\iunit 2\pi d_-\mb{ax+b}}}x^{-s-1}\on{d}x=-\sum_{k=1}^\infty F_{s,k}\mb{d_-a,d_-b,\mf 14,d_+}+\iunit F_{s,k}\mb{d_-a,d_-b,0,d_+}\\
			&=-\sum_{k=1}^\infty \mf{e^{\iunit 2\pi d_-bk}}{sd_+^{-s}k}{}_1F_1\mb{-s;1-s;\iunit 2\pi d_-ad_+k}+g_{s,k}\mb{d_-a,d_-b,\mf 14}+\iunit g_{s,k}\mb{d_-a,d_-b, 0}.
		}		
		Recalling the definition of $\mbs{\cdot}$ in \eq*{fractpart}, we factorize, using $\Re\mb{ 1-e^{\iunit 2\pi \Re z}}\ge 0$,
		\be*{align*}
		{
			\mbs{z}-\mf 12-\iunit \Im z&=\be*{cases}{\mf\iunit{2\pi}\log\mb{\mf{1-e^{-\iunit 2\pi \Re z}}{1-e^{\iunit 2\pi \Re z}}} & \Re z\in\MB{R}\setminus\MB{Z}\\-\mf 12 & \Re z\in \MB{Z}}\\
			&=\be*{cases}{-\mf\iunit{2\pi}\mb{\log\mb{1-e^{\iunit 2\pi \Re z}}-\log\mb{1-e^{-\iunit 2\pi \Re z}}} & \Re z\in\MB{R}\setminus\MB{Z}\\-\mf 12 & \Re z\in \MB{Z}}.
		}
		Using the Taylor series for $\log\mb{1-\cdot}$, we find now for almost all on $\MB{R}$ the Fourier series as stated in \eq*{fourierseries} (which is the sum of two convergent Dirichlet series). Hence we find, using \lem{mellintransform:I} and \lem{mellintransform:II}, assuming $-1<\Re s <0$,
		\be*{align*}
		{
			\pi\mb{\int_0^{d_+}-\int_0^d-\int_{d_+}^\infty}&\mb{\mf 12-\mbs{x\Re a+\Re b}}x^{-s-1}\on{d}x=\\
			\sum_{k=1}^\infty f_{s,k}&\mb{\Re a,\Re b,0,d}-f_{s,k}\mb{\Re a,\Re b,0,d_+}-F_{s,k}\mb{\Re a,\Re b,0,d_+}=\\
			\sum_{k=1}^\infty& f_{s,k}\mb{\Re a,\Re b,0,d}-2f_{s,k}\mb{\Re a,\Re b,0,d_+}-g_{s,k}\mb{\Re a,\Re b,0}.
		}
		Analogously, we compute
		\be*{align*}
		{
			\int_d^\infty\mb{\mf 12-\mbs{ax+b}}x^{-s-1}\on{d}x=\sum_{k=1}^\infty f_{s,k}\mb{\Re a,\Re b,0,d}+g_{s,k}\mb{\Re a,\Re b,0}.
		}
		To see the claim on $\Im a=\Im b=0$ note that then $d_+=d$ and $d_-=1$. Then the sum simplifies to (recalling trigonometric angle sum identities for $\sin$ and $\cos$)
		\be*{align*}
		{
			I_s\mb{a,b,c,d}=&\cos\mb{2\pi c}\sum_{k=1}^\infty f_{s,k}\mb{a,b,0,d}+g_{s,k}\mb{a,b,0}\\
			&+\sin\mb{2\pi c}\sum_{k=1}^\infty f_{s,k}\mb{a,b,\mf 14,d}+g_{s,k}\mb{a,b,\mf 14}\\
			=&\sum_{k=1}^\infty f_{s,k}\mb{a,b,c,d}+g_{s,k}\mb{a,b,c}.
		}
		To see the claim on holomorphicity, we know that all involved (or sums of it) $f_{s,k},F_{s,k}$, $k\in\MB{N}$ are holomorphic on the required domains, using \lem{mellintransform:II} and the bound in \eq*{mellintransform:I:bound:II} (meromorphic functions bounded at all possible poles are holomorphic) from \lem{mellintransform:I}. Then this claim follows by the uniform summable bounds on the $f_{s,k}$ and $F_{s,k}$ (or sums of it) in those required variables, invoking Weierstrass M-Test, given in \tm[V]{convergence}-\tm*[VI]{convergence}, \tm[VIII]{convergence}-\tm*[IX]{convergence}.	We now show the claimed bounds in \tm[V]{convergence}-\tm*[IX]{convergence}. We write
		\be*{align*}
		{
			2\iunit F_{s,k}\mb{a,b,c,d}=&F_{s,k}\mb{a,b,c+\mf 14,d}+\iunit F_{s,k}\mb{a,b,c,d}\\
			&-\mb{F_{s,k}\mb{-a,-b,-c+\mf14,d}+\iunit F_{s,k}\mb{-a,-b,-c,d}}\\
			=&e^{\iunit 2\pi \mb{bk+c}}\mb{F_{s,k}\mb{a,0,\mf 14,d}+\iunit F_{s,k}\mb{a,0,0,d}}\\
			&-e^{-\iunit 2\pi \mb{bk+c}}\mb{F_{s,k}\mb{-a,0,\mf 14,d}+\iunit F_{s,k}\mb{-a,0,0,d}}.
		}
		We fix $s\in\MB{C}$ and let $n\in\MB{N}$ be minimal such that $\Re s+n>0$. Then we invoke the recurrence relation together with the given bound from \lem{mellintransform:II:recurrence}, assuming $\Im a=\Im b=0$, to bound
		\be*{align*}
		{
			\mbm{F_{s,k}\mb{a,b,c,d}}\le& \cosh\mb{2\pi \Im\mb{adk+bk+c}}\\
			&\times{\sum_{m=0}^{n-1}\mf{\mbm{\mb{s+1}_m}}{\mb{2\pi\mbm{a}kd}^{m+1}kd^{\Re s}}+\mf{\mbm{\mb{s+1}_n}}{\mb{2\pi \mbm{a}k}^nk}\mf{d^{-\Re s-n}}{\Re s+n}}.
		}
		Hence on subsets of $\MB{C}$ such that $\Re s>-\epsilon$ and $\Im s\le \epsilon$ for every $\epsilon>0$ this gives a uniform bound in $s$ summable in $k$. 
		Conversely, we make us of the identity in \eq*{expidentity},
		\be*{align*}
		{
			&2\iunit\mb{f_{s,k}\mb{a,b,c,d_1}-f_{s,k}\mb{a,b,c,d_2}}=\\
			&\mf{e^{\iunit 2\pi \mb{bk+c}}}{sk}\mb{\mf{{}_1F_1\mb{-s,1-s,\iunit 2\pi ad_1k}}{d_1^s}-\mf{{}_1F_1\mb{-s,1-s,\iunit 2\pi ad_2k}}{d_2^s}}\\
			&-\mf{e^{-\iunit 2\pi \mb{bk+c}}}{sk}\mb{\mf{{}_1F_1\mb{-s,1-s,-\iunit 2\pi ad_1k}}{d_1^s}-\mf{{}_1F_1\mb{-s,1-s,-\iunit 2\pi ad_2k}}{d_2^s}},
		}
		and we bound the above \rhs, using the bound in \eq*{mellintransform:I:bound:II} from \lem{mellintransform:I}, with the analogous assumptions on $s,a,b,d_1,d_2$ as before.
		We now show the bounds claimed in \tm[V]{convergence}-\tm*[IX]{convergence}. The bound in \tm[IX]{convergence} follows directly by summation in $k$ of the identity from \lem{mellintransform:II:recurrence}, and using the bound for $\Re s+n>0$, taking $n$ large enough as before. The bound in \tm[VIII]{convergence} follows by summation in $k$ of the bound in \eq*{mellintransform:I:bound:II} from \lem{mellintransform:I:recurrence} in which we take $a\rightarrow  ad_-$, $d_1\rightarrow d$ and $d_2\rightarrow d_+$, assuming $\Im a \neq 0$ and $d_-\mb{d_+\Im a+\Im b}=0$,		
		\be*{align*}
		{
			&\sum_{k=1}^\infty \mbm{\mf {e^{-\iunit 2\pi d_-bk}}{sk}\mb{\mf {{}_1F_1\mb{-s;1-s;-\iunit 2\pi d_-adk}}{d^{1+s}}-\mf {{}_1F_1\mb{-s;1-s;-\iunit 2\pi d_-ad_+k}}{d_+^{1+s}}}}\le\\
			&\mf {L_2\mb{e^{2\pi d_-\mb{d_+\Im a+\Im b}}}}{2\pi\mbm{a}}\mb{\sup_{k\in\MB{N}}\mbm{\mf 1{d_+^{ s+1}}-\mf {e^{\iunit 2\pi d_-a\mb{d_+-d}k}}{d^{s+1}}}+\mf {\mbm{s+1}}{\mbm{1+\Re s}}\mbm{d^{-\Re s-1}-d_+^{-\Re s-1}}}.
		}
		Note that $L_2\mb{e^{2\pi d_-\mb{d_+\Im a+\Im b}}}=\zeta\mb{2}$ because of the assumption $d_-\mb{d_+\Im a+\Im b}=0$. In the case $d_-\mb{d_+\Im a+\Im b}\neq 0$ then $d_+=d$ and the above estimate is trivial.\\
		The bounds in \tm[VII]{convergence} follow by summation in $k$ of the recurrence relation from \lem{mellintransform:I:recurrence}, in which we take $n=1$, together with the bound in \eq*{mellintransform:I:bound:I}. Hence, assuming $d_-=1$ and $d_+>d$, we bound for every $d<x< d_+$
		\be{align}{eq:sum1F1bound}
		{
			&x\sum_{k=1}^\infty \mbm{\mf {e^{-\iunit 2\pi bk}}{sk} {{}_1F_1\mb{-s;1-s;-\iunit 2\pi axk}}}\le \mf{\on{Li}_2\mb{e^{2\pi \Im\mb{ad+b}}}}{2\pi\mbm{a}}+\mf{\on{Li}_2\mb{e^{ 2\pi k \Im b}}}{2\pi\mbm{a}}\\
			&+\mbm{s+1}\mb{\mf{\mb{1-e^x}}{\Re s}\mf{\on{Li}_{1-\Re s}\mb{e^{2\pi k \Im b}}}{\mb{2\pi\mbm{a}}^{-\Re s}}+2\mbm{\mf{\Phi_2\mb{e^{2\pi\Im b},k_0}}{2\pi\mbm{a}\mb{1+\Re s}}-\mf{\Phi_{1-\Re s}\mb{e^{2\pi\Im b},k_0}}{\mb{2\pi \mbm{a}}^{-\Re s}\mb{1+\Re s}}}}\nonumber\\
			&+2\mbm{s+1}\mbm{\mf{\Phi_2\mb{e^{2\pi\Im b},1}-\Phi_2\mb{e^{2\pi\Im b},k_0}}{2\pi\mbm{a}\mb{1+\Re s}}-\mf{\Phi_{1-\Re s}\mb{e^{2\pi\Im b},k_0}-\Phi_{1-\Re s}\mb{e^{2\pi\Im b},1}}{\mb{2\pi \mbm{a}}^{-\Re s}\mb{1+\Re s}}}\nonumber,
		}
		with the minimal choice (or $k=1$ if it does not exist) for $k_0\in\MB{N}$ such that
		\[\on{sgn}\mb{1-\mb{2\pi\mbm{a}\mb{k_0-1}}^{1+\Re s}}\neq \on{sgn}\mb{1-\mb{2\pi\mbm{a}k_0}^{1+\Re s}}.\]
		In addition, note that
		\[\Im b=\Im b+x\Im a-x\Im a\le 0,\]
		since $d_-=1$ and $d<x<d_+$. The case $d_-=-1$ and $\Im a=0$ is analogous, changing the appropriate signs for $a,b$.
		The bounds in \tm[V]{convergence}-\tm*[VI]{convergence} are deduced from the bound in \eq*{sum1F1bound}, recalling $f_{s,k}$ defined in \eq*{fk} and the identity in \eq*{expidentity}, noting
		\be*{align*}
		{
			&2\iunit kd^ssf_{s,k,1}\mb{\Re a,\Re b,c,d}=2\iunit kd^ss\cos\mb{2\pi \mb{\Re bk+c}}f_{s,k}\mb{\Re a,0,0,d}\\
			&=\cos\mb{2\pi \mb{\Re  bk+c}}\mb{{}_1F_1\mb{-s;1-s;\iunit 2\pi d\Re a }-{}_1F_1\mb{-s;1-s;-\iunit 2\pi d\Re a }},\\
			&2kd^ssf_{s,k,2}\mb{\Re a,\Re b,c,d}=2kd^ss\sin\mb{2\pi \mb{\Re bk+c}}f_{s,k}\mb{\Re a,0,\mf 14,d}\\
			&=\sin\mb{2\pi \mb{\Re  bk+c}}\mb{{}_1F_1\mb{-s;1-s;\iunit 2\pi d\Re a }+{}_1F_1\mb{-s;1-s;-\iunit 2\pi d\Re a }}.
		}
	}	
	
	We now prepare the proof of \tm{maintheorem}.
	
	Setting
	\[\mbfl{z}\coloneqq z-\mbs{z},\]
	we have
	\[\mbfl{z}= \mbfl{\Re z}.\]
	
	For every $s\in\MB{C}$, $a>0$, $b\in\MB{R}$, $d>0$ we set
	\be{align}{eq:remainder}
	{
		R_{a,b,d}\mb{s}&\coloneqq s\mb{\int_1^{\mf{1-\mbs{b}}a}-\int_1^d}\mbfl{ax+\mbs{b}}x^{-s-1}\on{d}x,
	}
	and in addition for every $\Re s>0$, $a\in\MB{R}\setminus\mbs{0}$, $c\in\MB{R}$ we set
	\be{align}{eq:Mellin}
	{
		M_{a,b,d}\mb{s}&\coloneqq s\int_d^\infty\mb{\mf 12-\mbs{ax+b}}x^{-s-1}\on{d}x.
	}
	
	\be{lemma}{l:functionidentity}
	{
		Let $s\in\MB{C}$ such that $\Re s>1$ and let $a\in\MB{R}\setminus\mbs{0}$, $b\in\MB{R}$, $d>0$. Then it holds
		\be*{align*}
		{
			M_{a,b,0,d}\mb{s}=\be*{cases}
			{
				a^s\zeta\mb{s,1-\mbs{b}}-\mb{\mf {ads}{1-s}+\mbs{b}-\mf 12}d^{-s}+R_{a,b,d}\mb{s},&a>0\\
				-\mb{-a}^s\zeta\mb{s,1-\mbs{-b}}+\mb{\mf {ads}{s-1}+\mbs{-b}-\mf 12}d^{-s}-R_{-a,-b,d}\mb{s},&a<0
			}.
		}
		In particular, the function
		\[a^s\zeta\mb{s,1-\mbs{b}}-\mb{\mf {ads}{1-s}+\mbs{b}-\mf 12}d^{-s}\] 
		extends holomorphically in $\Re s>0$ with the extension given by $M_{a,b,0,d}\mb{s}$.
	}
	\be*{proof}
	{
		Recalling the definition of $M_{a,b,d}$ in \eq*{Mellin}, we calculate, using first $\Re s>1$ and $a>0$,
		\be*{align*}
		{
			&s\int_d^\infty \mb{\mf 12-\mbs{ax+b}}x^{-s-1}\on{d}x=s\int_d^\infty \mb{\mbfl{ax+b}+\mf12-ax-b}x^{-s-1}\on{d}x\\
			&=s\int_d^\infty \mb{\mbfl{ax+\mbs{b}}+\mf12-ax-\mbs{b}}x^{-s-1}\on{d}x\\
			&=a^s\zeta\mb{s,1-\mbs{b}}+\mb{\mf 12-\mbs{b}}d^{-s}+ad^{1-s}\mf s{s-1}+s\mb{\int_d^\infty-\int_{\mf {1-\mbs{b}}a}^\infty}\mf{\mbfl{ax+\mbs{b}}}{x^{s+1}}\on{d}x\\
			&=a^s\zeta\mb{s,1-\mbs{b}}+\mb{\mf 12-\mbs{b}}d^{-s}+ad^{1-s}\mf s{s-1}+R_{a,b,d}\mb{s}.
		}
		Since the \lhs extends holomorphically on $\Re s>0$ and since $R_{a,b,d}$ is holomorphic on $\MB{C}$, we conclude for $a>0$. To see the result for $a<0$, we note for almost all $x\in\MB{R}$
		\[\mf 12-\mbs{ax+b}=-\mb{\mf 12-\mbs{-ax-b}}.\]
	}
	
	\be*[Proof of \tm{maintheorem}]{proof}
	{
		The result on the principal identity follows, using \lem{functionidentity} together with \tm{convergence} in which we take $a\in\MB{R}\setminus\mbs{0}$, $b\in\MB{R}$, $c=0$, $d>0$, noting
		\[I_s\mb{a,b,0,d}=\sum_{k=1}^\infty f_{s,k}\mb{a,b,0,d}+g_{s,k}\mb{a,b,0}.\]
		To see the statement on the two identities for $\Re s<0$ and $\Re s <1$, we use the bounds in \tm[V]{convergence}-\tm*[VI]{convergence}, noting
		\[\lim_{d\to 0}\sum_{k=1}^\infty f_{s,k}\mb{a,b,0,d}=0,\]
		if $\Re s<-1$ and if $d\le \mf {1-\mbs{b}}a$ it holds
		\[R\mb{a,b,d}\equiv 0.\]
		Comparing both sides in the principal identity, we find, assuming $a>0$,
		\[a^s\mf \pi s \zeta\mb{s,1-\mbs{b}}=\sum_{k=1}^\infty k^{s-1}\mb{2\pi}^s\Gamma\mb{-s}\mf 1{2\iunit}\mb{\mb{-\iunit a}^se^{\iunit 2\pi{bk}}-\mb{\iunit a}^se^{-\iunit 2\pi{bk}}}, \]
		and its \rhs converges in $\Re s<0$ and in $\Re s<1$ if $b=0$ and $b\neq 0$, respectively.
	}
	
	We now prepare the proof of \tm{maintheorem:B}.

	\be{lemma}{l:mellintransform:III:recurrence}
	{
		Let $s\in\MB{C}\in\MB{N}$, $a,b\in\MB{C}$ and let $d>0$. Then for every $n\in\MB{N}$ and for every $k\in\MB{Z}$ such that $ad+b\neq -k\neq b$ it holds
		\be*{align*}
		{
			\mf{ad}{1-s}\mf{{}_2F_1\mb{1,1-s;2-s;-\mf{ad}{b+k}}}{b+k}=&\sum_{m=1}^{n-1}\mf{\mb{ad}^{m}\mb{m-1}!}{\mb{1-s}_{m}\mb{ad+b+k}^{m}}\\
			&+\mf{\mb{ad}^{n}{}_2F_1\mb{n,-s+n;1-s+n;-\mf{ad}{b+k}}}{\mb{1-s}_{n}\mb{b+k}^{n}}.
		}
		Moreover, for every $n\in\MB{N}$ such that $n-\Re s>0$, for every $\epsilon>0$ and for every $k\in\MB{Z}$ such that $\mbm{\mf{ad}{b+k}}<\epsilon$ it holds
		\be{align}{eq:bound:mellin:III}
		{
			\mbm{\mb{ad}^{n}{}_2F_1\mb{n,-s+n;1-s+n;-\mf{ad}{b+k}}}\le \mbm{\mf{s-n}{\Re s-n}}\mbm{\mf {ad}{b+k}}^n\mbm{1-\epsilon}^{-n}.
		}
		Moreover, it holds
		\be*{align*}
		{
			\mf{ad}{1-s}\mf{{}_2F_1\mb{1,1-s;2-s;-\mf{ad}{b+k}}}{b+k}=\mf{ad}{1-s}\mf{{}_2F_1\mb{1,1;2-s;\mb{1+\mf{b+k}{ad}}^{-1}}}{ad+b+k}.		
		}
	}
	\be*{proof}
	{
		Using Euler's formula \cite[p.57, Equation (1)]{Luke_1969}
		\[
		{}_2F_1\mb{n,-s+n;1-s+n;-\mf{ad}{b+k}}=\mbm{\int_0^1\mf{x^{-s-1}}{\mb{adx+b+k}^n}\on{d}x}
		\]
		for the first, noting that their $b$ is a natural number in our case hence it is enough to assume $ad+b\neq -k\neq b$, and integration by parts for the second identity, we compute
		\be*{align*}
		{
			\mf{d^{1-s}}{1-s}\mf{{}_2F_1\mb{1,1-s;2-s;-\mf{ad}{b+k}}}{b+k}&=\int_0^d\mf{x^{-s}}{ax+b+k}\on{d}x\\
			=\mf{\mb{ad+b+k}^{-1}d^{1-s}}{1-s}&+\mf{a}{1-s}\int_0^d\mf{x^{1-s}}{\mb{ax+b+k}^2}\on{d}x,
		}
		and, repeating integration by parts $n-1$ times, we find
		\be*{align*}
		{
			{}_2F_1\mb{1,1-s;2-s;-\mf{ad}{b+k}}=&\sum_{m=1}^{n-1}\mf{\mb{b+k}\mb{ad}^{m-1}\mb{m-1}!}{\mb{2-s}_{m-1}\mb{ad+b+k}^{m}}\\
			&+\mf{\mb{ad}^{n-1}{}_2F_1\mb{n,-s+n;1-s+n;-\mf{ad}{b+k}}}{\mb{2-s}_{n-1}\mb{b+k}^{n-1}}.
		}
		The claimed upper bound follows, estimating, using $n-\Re s>0$ for the second inequality,
		\be*{align*}
		{
			&\mbm{\mb{b+k}^{-n}{}_2F_1\mb{n,-s+n;1-s+n;-\mf{ad}{b+k}}}=\mbm{\mb{s-n}\int_0^1\mf{x^{n-s-1}}{\mb{adx+b+k}^n}\on{d}x}\\
			&\mbm{\mf{s-n}{\Re s-n}}\mbm{{b+k}}^{-n}\sup_{0<x<1}\mbm{1+x\mf{ad}{b+k}}^{-n}\le \mbm{\mf{s-n}{\Re s-n}}\mbm{{b+k}}^{-n}\mb{1-\epsilon}^{-n}.
		}
	}
	
	\be{lemma}{l:mellintransform:III}
	{
		Let $s\in\MB{C}$ such that $\Re s<0$, let $a,b\in\MB{C}$ and let $d>0$. Let $n_1,n_2\in\MB{Z}$ be the smallest and largest integer such that 
		\[0<\mf{n_1-\Re b\on{sgn}\Re a}{\mbm{\Re a}}\le \mf{n_2-\Re b\on{sgn}\Re a}{\mbm{\Re a}}< d,\]
		respectively, and otherwise set $n_1\coloneqq 1$, $n_2\coloneqq 0$. Let $n_0\in\MB{Z}$ be such that
		\[\on{sgn}\mb{\Im a\mf{n_0-\Re b\on{sgn}\Re a}{\mbm{\Re a}}+\Im b}\neq \on{sgn}\mb{\Im a\mb{\mf{n_0+1-\Re b\on{sgn}\Re a}{\mbm{\Re a}}}+\Im b},\]
		and if no such $n_0$ exists or $n_0<n_1$ set $n_0\coloneq n_1-1$ and otherwise set $n_1=n_2$. Set
		\[\omega_-\coloneqq1_{|\MB{R}\setminus\mbs{0}}\mb{\Re a}\be*{cases}{\on{sgn}\Re a\on{sgn}\Im a&\ft{if }\Im a\neq 0\\
			\on{sgn}\Re a\on{sgn}\Im b&\ft{if }\Im a=0\neq \Im b\\
			-1&\ft{if }\Im a=\Im b=0
		}.\]
		If either $ad+b\not\in\MB{Z}$ or $ad+b\in\MB{Z}$ and $a\neq 0$ then it holds
		\be*{align*}
		{
			\int_0^d&\log\mb{2\sin\mb{\pi\mb{ax+b}}}^2x^{-s-1}\on{d}x=-1_{|\MB{R}\setminus\MB{Z}}\mb{\Re\mb{ad+b}}\mf {\log\mb{\sin\mb{\pi\mb{ad+b}}}^2}{sd^s}\\
			&-\mf{1}{sd^s}1_{|\MB{Z}+\iunit\MB{R}\setminus\mbs{0}}\mb{ad+b}\mb{\log{\mb{\sin\mb{\pi\mb{ad+b}}}^2}-1_{|\mbs{1}}\mb{\omega_-}\iunit 2\pi}\\
			&-\mf{2\log 2}{sd^s}+\mf {2}{sd^s}1_{|\MB{Z}}\mb{ad+b}\mb{\gamma+\psi\mb{1-s}-\mf 12\log\mb{\pi ad}^2}-1_{|\MB{Z}}\mb{b}\mf {2}{s^2d^s}\\
			&-\omega_-\iunit \mf{2\pi}s\mbm{\Re a}^s\mb{\zeta\mb{s,n_1-\Re b\on{sgn}\Re a}-\zeta\mb{s,n_0+1-\Re b\on{sgn}\Re a}}\\
			&+\omega_-\iunit \mf{2\pi}s\mbm{\Re a}^s\mb{\zeta\mb{s,n_0+1-\Re b\on{sgn}\Re a}-\zeta\mb{s,n_2+1-\Re b\on{sgn}\Re a}}\\
			&+\mf{2a}{d^{s-1}}\sum_{\substack{k\in\MB{Z}\\k\not\in\mbs{b,ad+b}}}\mf{{}_2F_1\mb{1,1-s;2-s;-\mf{ad}{b-k}}}{s\mb{1-s}\mb{b-k}},
		}
		where $\gamma$ denotes the Euler-Mascheroni constant and $\psi$ the Digamma function. Moreover, the above infinite sum converges to a meromorphic function in $s\in\MB{C}$ with poles in $s\in\MB{Z}_{\ge 2}$. The convergence is uniform in $s$ for every $\epsilon>0$ such that $\Re s <\epsilon$ and $s$ bounded away from $\MB{N}_{\le \epsilon}$, uniform in $b$ uniformly bounded away from $\MB{Z}$ and uniform in $a$ and $d$ with $ad+b$ bounded away from $\MB{Z}$ and $ad$ bounded.
	}
	\be*{proof}
	{
		We note first for every $z\in\MB{C}\setminus\MB{Z}$
		\[\log\mb{\sin \mb{\pi z}}^2=\lim_{\epsilon\to 0}\log\mb{\sin \mb{\pi z\pm\iunit\epsilon}}^2.\]
		A point $x\in\MB{R}$ of discontinuity \wrt the principal branch for $\log$ for
		\[\Im\log\mb{\sin\mb{\pi\mb{ax+b\pm\iunit\epsilon}}}^2\] 
		is given whenever $\sin\mb{\pi\mb{ax+b\pm\iunit\epsilon}}\in\iunit \MB{R}$ which is equivalent to $\Re\mb{ax+b}\in\MB{Z}$, assuming $\Re a\neq 0$. Recalling the assumptions for $n_1,n_2$, we set
		for every $m\in\MB{N}_{\le {n_2-n_1+1}}$
		\[x_0\coloneq 0,\quad x_{m}\coloneqq \mf{n_1+m-1-\Re b\on{sgn}\Re a}{\mbm{\Re a}},\quad x_{n_2-n_1+2}\coloneqq d.\]
		We set for every $\epsilon>0$ and for every $1\le m\le n_2-n_1+1$
		\[\epsilon_a=\epsilon_a\mb{\epsilon}\coloneqq \on{sgn}\mb{\Re a}\epsilon\ftqq y_m\coloneqq \lim_{\epsilon\to 0}\lim_{\delta\to 0}\Im \log\mb{\sin\mb{\pi\mb{a\mb{x_m+\delta}+b-\iunit \epsilon_a}}}^2.\]
		Then, assuming $\Re s <0$, using absolute convergence of the integral uniformly bounded in $\epsilon$ to exchange the limit, we write
		\be*{align*}
		{
			\int_0^d\mf{\log\mb{2\sin\mb{\pi\mb{ax+b}}}^2}{x^{s+1}}\on{d}x&=\lim_{\epsilon\to 0}\int_0^d\mf{\log\mb{2\sin\mb{\pi\mb{ax+b-\iunit\epsilon_a}}}^2}{x^{s+1}}\on{d}x\\
			=\lim_{\epsilon\to 0}\sum_{m=0}^{n_2-n_1+1}\int_{x_m}^{x_{m+1}}&\mf{\log\mb{2\sin\mb{\pi\mb{ax\on{sgn}\Re a+b\on{sgn}\Re a-\iunit\epsilon}}}^2}{x^{s+1}}\on{d}x.
		}
		We calculate, recalling $\Re s<0$, assuming $ad+b,b\not\in \MB{Z}$ and for now $\Re a\ge 0$ since the result is invariant under the simultaneous change $a\to a\on{sgn}\Re a$ and $b\to b\on{sgn}\Re a$, integrating by parts for the first equality, using absolute convergence for the second equality as the sum is understood in summing in symmetric sets in $\MB{Z}$, and using Euler's formula \cite[p.57, Equation (1)]{Luke_1969} for the third equality,
		\be*{align*}
		{
			&\int_{x_m}^{x_{m+1}}\mf{\log\mb{2\sin\mb{\pi\mb{ax+b-\iunit\epsilon}}}^2}{x^{s+1}}\on{d}x+\\
			&\lim_{\delta\to 0}\mf {\log\mb{2\sin\mb{\pi\mb{a\mb{x_{m+1}-\delta}+b-\iunit\epsilon}}}^2}{s\mb{x_{m+1}-\delta}^s}-\mf {\log\mb{2\sin\mb{\pi\mb{a\mb{x_m+\delta}+b-\iunit\epsilon}}}^2}{s\mb{x_m+\delta}^s}\\
			&=\int_{x_m}^{x_{m+1}}\mf{2\pi as^{-1}x^{-s}}{\tan\mb{\pi\mb{ax+b-\iunit\epsilon}}}\on{d}x=\mf {2a}s\sum_{k\in\MB{Z}}\int_{x_m}^{x_{m+1}}\mf{x^{-s}}{ax+b-\iunit\epsilon+k}\on{d}x\\
			&=2a\sum_{k\in\MB{Z}}\mf{x_{m+1}^{1-s}{}_2F_1\mb{1,1-s;2-s;-\mf{ax_{m+1}}{b-\iunit\epsilon+k}}-x_{m}^{1-s}{}_2F_1\mb{1,1-s;2-s;-\mf{ax_{m}}{b-\iunit\epsilon+k}}}{{s\mb{1-s}}\mb{b-\iunit\epsilon+k}},
		}
		where the term for $x_0$ here and in the rest of the proof is understood as $0$ since $\Re s<0$. In the case $b\in\MB{Z}$ then the corresponding integral terms for $k=-b$ simplify to
		\be{align}{eq:mod:I}
		{
			{2a}s\int_{x_m}^{x_{m+1}}\mf{x^{-s}}{ax+b+k}\on{d}x=\mf 2s\int_{x_m}^{x_{m+1}}x^{-s-1}\on{d}x=\mf {2\mb{x_{m}^{-s}-x_{m+1}^{-s}}}{s^2}.
		}
		In the case $ad+b\in\MB{Z}$ we write the corresponding integral
		\be{align}{eq:mod:II}
		{
			&\int_0^d\mf{\log\mb{2\sin\mb{\pi\mb{ax+b}}}^2}{x^{s+1}}\on{d}x=			\int_0^d\mf{\log\mb{2\sin\mb{\pi{a\mb{x-d}}}}^2}{x^{s+1}}\on{d}x\nonumber\\
			&=\int_0^d\mf{\log\mb{\mf{2\sin\mb{\pi a\mb{x-d}}}{\pi \mbm{a}\mb{x-d}}}^2}{x^{s+1}}\on{d}x-\mf {\log\mb{\pi \mbm{a}d}^2}{sd^s}+2\int_0^d\mf{\log\mb{1-\mf xd}}{x^{s+1}}\on{d}x,
		}
		which is valid if $a\neq 0$. For the first integral we proceed analogously as in the case $ad+b\not\in \MB{Z}$, noting that the corresponding points of discontinuity are not altered, and the second integral we compute directly
		\be{align}{eq:mod:III}
		{
			\int_0^d\log\mb{1-\mf xd}{x^{-s-1}}\on{d}x=-d^{-s}\sum_{k=1}^\infty \mf 1{k\mb{k-s}}=\mf{\psi\mb{1-s}+\gamma}{sd^s},
		}
		where the last equality follows from the Weiserstrass factorization for the reciprocal Gamma function $1\slash\Gamma$ and the definition for the Digamma function $\psi$ (e.g. see \cite[pp. 47]{Remmert_2007}. To see the statement on convergence and meromorphicity we use the recurrence relation from \lem{mellintransform:III:recurrence} in which we take $n>\Re s$ together with the bound in \eq*{bound:mellin:III} which is valid due to the assumptions on $s,a,b$ and $d$. This allows us to find meromorphicity in $s\in\MB{C}$ of the summands for every sufficiently large $k$ in question with simple poles in $s\in\MB{Z}_{\ge 0}$. (Note that there are only finitely many summands with smaller values of $k$ due to the assumed bounds for $a,b,d$.) Invoking the Weiserstrass M-Test, the convergence is uniform accordingly and hence the sum converges to a meromorphic function with prescribes simple poles.
		Now we sum in $m$, and by telescoping, assuming again $ad+b,b\not\in \MB{Z}$ (the other cases are treated accordingly and pose only minor modifications which follow directly from our above computations in \eq*{mod:I}-\eq*{mod:III}), we arrive at
		\be*{align*}
		{
			\int_0^d\mf{\log\mb{2\sin\mb{\pi\mb{ax+b}}}^2}{x^{s+1}}\on{d}x&+\mf {\log\mb{2\sin\mb{\pi\mb{ad+b}}}^2}{sd^s}=\\
			\iunit 2\sum_{m=1}^{n_2-n_1+1}\lim_{\epsilon\to 0}\lim_{\delta\to 0}&\Im\mf {\log\mb{\sin\mb{\pi\mb{a\mb{x_m+\delta}+b-\iunit\epsilon_a}}}^2}{sx_m^s}\\
			&+\mf {2a}{d^{s-1}}\lim_{\epsilon\to 0}\sum_{k\in\MB{Z}}\mf{{}_2F_1\mb{1,1-s;2-s;-\mf{ad}{b-\iunit\epsilon_a+k}}}{s\mb{1-s}\mb{b-\iunit\epsilon_a +k}}\\
			=&\sum_{m=1}^{n_2-n_1+1}\iunit2\mf{y_m}{sx_m^s}+\mf {2a}{d^{s-1}}\sum_{k\in\MB{Z}}\mf{{}_2F_1\mb{1,1-s;2-s;-\mf{ad}{b+k}}}{s\mb{1-s}\mb{b+k}},		
		}
		where the last equality follows by continuity for ${}_2F_1\mb{1,1-s,2-s,\cdot}$ on $\MB{C}\setminus\mbs{1}$ in the directions belonging to the closed upper hemi sphere. Now $n_2-n_1+1=0$ if $\Re a=0$. Assuming $\Re a\neq 0$, then $y_m$ is just the argument at the discontinuity from the right, meaning it is either $\pi$ or $-\pi$ depending on the orientation at the intersection and the location on the imaginary axis of $\sin\mb{a\mb{x_m+\delta}+b-\iunit \epsilon_a}$. Hence it holds for every $\epsilon>0$ small enough
		\[y_m=\pi\on{sgn}\mb{\lim_{\delta\to 0}\Im \log\mb{\sin\mb{a\mb{x_m+\delta}+b-\iunit \epsilon_a}}^2}=\pi\omega_-.\]
		Then we calculate, recalling $n=0$ if $\Re a=0$, 
		\be*{align*}
		{
			\sum_{m=1}^\infty\mf{y_m}{sx_m^s}=&\omega_-\mf{\pi}{s}\mb{\sum_{m=n_0+1}^{n_2}-\sum_{m=n_1}^{n_0}}\mb{\mf{m-\Re b\on{sgn}\Re a}{\mbm{\Re a}}}^{-s}\\
			=&\omega_-\mf{\pi}{s}\mbm{\Re a}^s\mb{\zeta\mb{s,n_0+1-\Re b\on{sgn}\Re a}-\zeta\mb{s,n_2+1-\Re b\on{sgn}\Re a}}\\
			&-\omega_-\mf{\pi}{s}\mbm{\Re a}^s\mb{\zeta\mb{s,n_1-\Re b\on{sgn}\Re a}-\zeta\mb{s,n_0+1-\Re b\on{sgn}\Re a}}.
		}
	}

	\be{lemma}{l:hyperidentity}
	{
		For every $s\in\MB{C}\setminus\MB{Z}$ and for every $z\in\MB{C}\setminus\mbs{0,-1}$ it holds
		\be{align}{eq:hyperidentity:II}
		{
			z\mf{{}_2F_1\mb{1,1-s;2-s,-z}}{1-s}+\mf{{}_2F_1\mb{1,s;1+s,-\mf 1z}}{s}=\mf\pi{\sin\mb{\pi s}}\be*{cases}
			{
				z^s&\ft{if }z\in\MB{C}\setminus\mbb{-1,0}\\
				\overline{\bar{z}^{\bar{s}}}&\ft{if }z\in\MB{C}\setminus\mblo{-\infty,-1}}.
		}
	}
	\be*{proof}
	{
		The claim in the case $z\in\MB{C}\setminus\mbb{-1,0}$ is the identity in \eq*{hyperidentity}. In particular, setting for every $x>1$ and every $\epsilon>0$
		\[z_\epsilon\coloneqq -x+\iunit\epsilon,\]
		the identity in \eq*{hyperidentity} with $z\rightarrow z_\epsilon$ holds in the limit $\epsilon\to 0$ by continuity of ${}_2F_1\mb{\dots;-z}$ and $z^s$ in $z$ in the direction of the closed upper and lower hemisphere, respectively. The term ${}_2F_1\mb{\dots;-\mf 1z}$ is holomorphic in $z$ since $x>1$. Now, consider the involution $z\mapsto \mf 1z$, which amounts to passing from $x>1$ to $0<x<1$.
		Then again in \eq*{hyperidentity} we have a term ${}_2F_1\mb{\dots;-z}$ and ${}_2F_1\mb{\dots;-\mf 1z}$ and by our above considerations, everything is continuous in the appropriate sense. However, on the \lhs we find
		\[\mb{z^{-1}}^s=\mb{\mf{-x-\iunit \epsilon}{x^2+\epsilon^2}}^s,\]
		which is not continuous in the limit in the direction we are taking here. Instead we can use
		\[\overline{\mb{\overline{\mb{z^{-1}}}}^{\bar s}},\]
		which takes the limit in the direction of the closed lower hemisphere. Since the limit on the \rhs is unique, we conclude.
	}
	
	\be{lemma}{l:mellintransform:IV}
	{
		Let $a,b\in\MB{C}$ and let $d>0$. Let $n\in\MB{Z}$ be the smallest integer such that 
		\[{n-\Re b\on{sgn}\Re a}{}> d\mbm{\Re a}.\]
		Let $n_0\in\MB{Z}$ be such that
		\[\on{sgn}\mb{\Im a\mb{\mf{n_0-\Re b\on{sgn}\Re a}{\mbm{\Re a}}}+\Im b}\neq \on{sgn}\mb{\Im a\mb{\mf{n_0+1-\Re b\on{sgn}\Re a}{\mbm{\Re a}}}+\Im b},\]
		and if no such $n_0$ exists or $n_0<n$ set $n_0\coloneqq n-1$. Set
		\[\omega_+\coloneqq1_{|\MB{R}\setminus\mbs{0}}\mb{\Re a}\be*{cases}{\on{sgn}\Re a\on{sgn}\Im a&\ft{if }\Im a\neq 0\\
			\on{sgn}\Re a\on{sgn}\Im b&\ft{if }\Im a=0\neq \Im b\\
			1&\ft{if }\Im a=\Im b=0
		}.\]
		$s\in\MB{C}$ such that $\Re s>1$ and if either $ad+b\not\in\MB{Z}$ or $ad+b\in\MB{Z}$ and $a\neq 0$ then it holds
		\be*{align*}
		{
			\int_d^\infty&\log\mb{2\sin\mb{\pi\mb{ax+b}}}^2x^{-s-1}\on{d}x=1_{|\MB{R}\setminus\MB{Z}}\mb{\Re\mb{ad+b}}\mf {\log\mb{\sin\mb{\pi\mb{ad+b}}}^2}{sd^s}\\
			&+\mf{1}{sd^s}1_{|\MB{Z}+\iunit\MB{R}\setminus\mbs{0}}\mb{ad+b}\mb{\log{\mb{\sin\mb{\pi\mb{ad+b}}}^2}-1_{|\mbs{-1}}\mb{\omega_+}\iunit 2\pi}\\
			&+\mf{2\log 2}{sd^s}+\mf {2}{sd^s}1_{|\MB{Z}}\mb{ad+b}\mb{\mf {1}{s}-{\psi\mb{1+s}-\gamma}+\mf 12\log\mb{\pi ad}^2}+1_{|\MB{Z}}\mb{b}\mf {2}{s^2d^s}\\
			&+\omega_+\iunit \mf{2\pi}{s}\mbm{\Re a}^s\mb{2\zeta\mb{s,n_0+1-\Re b\on{sgn}\Re a}-\zeta\mb{s,n-\Re b\on{sgn}\Re a}}\\
			&-\mf{2}{d^{s}}\sum_{\substack{k\in\MB{Z}\\k\not\in\mbs{b,ad+b}}}\mf{ad}{b-k}\mf{{}_2F_1\mb{1,1-s;2-s;-\mf{ad}{b-k}}}{s\mb{1-s}}\\
			&\quad\quad-\mf \pi{s\sin\mb{\pi s}}\mb{1_{|\MB{C}\setminus\MB{R}_{<-1}}\mb{\mf{b-k}{ad}}\mb{\mf{ad}{b-k}}^s+1_{|\MB{R}_{<-1}}\mb{\mf{b-k}{ad}}\overline{\mb{\overline{\mf{ad}{b-k}}}^{\bar s}}}.
		}
		where $\gamma$ denotes the Euler-Mascheroni constant and $\psi$ the Digamma function.
	}
	\be*{proof}
	{
		We proceed analogously to the proof for \lem{mellintransform:III}. In what follows we recalculate only the key steps and we have to make sense of the possible appearing infinite series because of infinitely many discontinuities. To this we set
		for every $m\in\MB{N}$
		\[x_{0}\coloneq d,\quad x_{m}\coloneqq \mf{n+m-\Re b\on{sgn}\Re a}{\mbm{\Re a}},\]
		and we set for every $\epsilon>0$ and for every $m\ge 1$
		\[\epsilon_a\coloneqq \epsilon\on{sgn}\Re a\ftqq y_m\coloneqq \lim_{\epsilon\to 0}\lim_{\delta\to 0}\Im \log\mb{\sin\mb{\pi\mb{a\mb{x_m+\delta}+b+\iunit \epsilon_a}}}^2.\]
		Again the exchange of the limit in $\epsilon$ with corresponding sums or integrals is justified by uniform convergence in $\epsilon$, using now $\Re s>0$. The sum over the functions ${}_1F_2$ enters, making use of
		\be*{align}
		{
			&\int_{x_m}^{x_{m+1}}\mf{\pi as^{-1}x^{-s}}{\tan\mb{\pi\mb{ax+b+\iunit\epsilon_a}}}\on{d}x=			\int_{x_{m+1}^{-1}}^{x_{m}^{-1}}\mf{\pi as^{-1}x^{s-2}}{\tan\mb{\pi\mb{\mf ax+b+\iunit\epsilon_a}}}\on{d}x\nonumber\\
			&=\sum_{k\in\MB{Z}}\int_0^1\mb{\mf{x_m^{-s}}{{1+x\mf{b+\iunit\epsilon_a+k}{ax_m}}}-\mf{x_{m+1}^{-s}}{{1+x\mf{b+\iunit\epsilon_a+k}{ax_{m+1}}}}}x^{s-1}\on{d}x\nonumber\\
			&=\sum_{k\in\MB{Z}}\mf{{}_2F_1\mb{1,s;1+s;-\mf{b+\iunit\epsilon_a+k}{ax_{m}}}}{sx_{m}^{s}}-\mf{{}_2F_1\mb{1,s;1+s;-\mf{b+\iunit\epsilon_a+k}{ax_{m+1}}}}{sx_{m+1}^{s}}.
		}
		We pass then from ${}_2F_1\mb{1,s;1+s;\cdot}$ to ${}_2F_1\mb{1,1-s;2-s;\cdot}$, using the identity in \eq*{hyperidentity:II} from \lem{hyperidentity}.	In the case $b\in\MB{Z}$ the corresponding integral terms for $k=-b$ simplify to
		\be{align}{eq:mod:IV}
		{
			\int_0^1\mb{\mf{x_m^{-s}}{{1+x\mf{b+k}{ax_m}}}-\mf{x_{m+1}^{-s}}{{1+x\mf{b+k}{ax_{m+1}}}}}x^{s-1}\on{d}x&=\int_0^1\mb{x_m^{-s}-x_{m+1}^{-s}}x^{s-1}\on{d}x\nonumber\\
			&=\mf {x_m^{-s}-x_{m+1}^{-s}}s.
		}
		In the case $ad+b\in\MB{Z}$ we write the corresponding integral
		\be{align}{eq:mod:V}
		{
			&\int_d^\infty\mf{\log\mb{2\sin\mb{\pi\mb{ax+b}}}^2}{x^{s+1}}\on{d}x=			\int_d^\infty\mf{\log\mb{2\sin\mb{\pi{a\mb{x-d}}}}^2}{x^{s+1}}\on{d}x\nonumber\\
			&=\int_d^\infty\mf{\log\mb{\mf{2\sin\mb{\pi a\mb{x-d}}}{\pi \mbm{a}\mb{x-d}}}^2}{x^{s+1}}\on{d}x+\mf {\log\mb{\pi \mbm{a}d}^2}{sd^s}+2\int_d^\infty\mf{\log\mb{\mf xd-1}}{x^{s+1}}\on{d}x,
		}
		which is valid if $a\neq 0$. For the first integral we proceed analogously as in the case $ad+b\not\in \MB{Z}$, noting that the corresponding points of discontinuity are not altered, and the second integral we compute directly
		\be{align}{eq:mod:VI}
		{
			&\int_d^\infty\log\mb{\mf xd-1}{x^{-s-1}}\on{d}x=d^{-s}\int_1^\infty\mb{\log x+\log\mb{1-\mf 1x}}{x^{-s-1}}\on{d}x\nonumber=\\
			&=\mf {d^{-s}}{s^2}-d^{-s}\sum_{k=1}^\infty \mf 1{k\mb{k+s}}=\mf {d^{-s}}{s^2}-\mf{\psi\mb{1+s}+\gamma}{sd^s}.
		}
		Hence we are left in dealing with the infinite sum over the discontinuities. At first, we find
		\[y_m=\pi\omega_+.\]
		Then we calculate, recalling $n=0$ if $\Re a=0$,
		\be*{align*}
		{
			\sum_{m=1}^\infty\mf{y_m}{sx_m^s}&=\omega_+\mf{\pi}{s}\mb{\sum_{m=n_0-n}^\infty-\sum_{m=1}^{n_0-n-1}}\mb{\mf{n+m-\Re b\on{sgn}\Re a}{\mbm{\Re a}}}^{-s}\\
			&=\omega_+\mf{\pi}{s}\mbm{\Re a}^s\mb{2\zeta\mb{s,n_0-\Re b\on{sgn}\Re a}-\zeta\mb{s,n+1-\Re b\on{sgn}\Re a}}.
		}
	}
	
	\be*[Proof of \tm{maintheorem:B}]{proof}
	{
		The result of the principal identity for $s\in\MB{C}$, $\Re s>1$, follows, using \tm{convergence} in which we take $a\in\MB{R}\setminus\mbs{0}$, $b\in\MB{R}$, $c=\mf 14$, $d>0$, noting
		\[I_s\mb{a,b,\mf 14,d}=\sum_{k=1}^\infty f_{s,k}\mb{a,b,\mf 14,d}+g_{s,k}\mb{a,b, \mf 14},\]
		together with the result of \lem{mellintransform:IV} \iww $a\in\MB{R}\setminus\mbs{0}$, $b\in\MB{R}$, $d>0$ to evaluate the \lhs.
		We recall that $n$ in \lem{mellintransform:IV} is the minimal integer such that $n>\mbm{a}d+b\on{sgn}{a}$ and we let $n_2\in\MB{Z}$ be maximal such that $n_2<b\on{sgn}{a}$.	We are left with the evaluation, assuming first $b\not\in\MB{Z}$,
		\be*{align*}
		{
			\sum_{\substack{k\in\MB{Z}\\k\neq\mbm{a}d+b\on{sgn}a}}1_{|\MB{R}\setminus\MB{R}_{<-1}}&\mb{\mf{b-k\on{sgn}a}{ad}}\mb{\mf 1{b\on{sgn} a-k}}^{s}\\
			&+1_{|\MB{R}_{<-1}}\mb{\mf{b-k\on{sgn}a}{ad}}\overline{\mb{\mf 1{b\on{sgn} a-k}}^{\bar s}}=\\
			\zeta\mb{s,b\on{sgn}a-n_2}&+e^{-\iunit\pi s}\zeta\mb{s,n-b\on{sgn}a}\\
			&+e^{\iunit\pi s}\mb{\zeta\mb{s,n_2+1-b\on{sgn}a}-\zeta\mb{s,n-b\on{sgn}a}}=\\
			\zeta\mb{s,1-\mbs{-b\on{sgn}a}}&+e^{\iunit\pi s}\zeta\mb{s,1-\mbs{b\on{sgn}a}}-2\iunit \sin\mb{\pi s}\zeta\mb{s,n-b\on{sgn}a}=\\
			\zeta\mb{s,1-\mbs{-b\on{sgn}a}}&+\cos\mb{\pi s}\zeta\mb{s,1-\mbs{b\on{sgn}a}}\\
			&+\iunit \sin\mb{\pi s}\mb{\zeta\mb{s,1-\mbs{b\on{sgn}a}}-2\zeta\mb{s,n-b\on{sgn}a}}.
		}
		If $b\in\MB{Z}$ we have to exclude the term for $k=b\on{sgn}a$. In particular, it holds $b\on{sgn}a-n_2=1=1-\mbs{-b\on{sgn}a}$ and we have to shift $n_2$ by $2$ instead of $1$ for the third term, meaning we replace in the first equality $n_2+1-b\on{sgn}a$ by $n_2+2-b\on{sgn}a=1=1-\mbs{b\on{sgn}a}$. Hence the above calculation is valid in the case $b\in\MB{Z}$ as well and we conclude.
		The statement about holomorphicity and convergence on the \rhs and the second identity in the case $\Re s<0$, $b\in\MB{Z}$ or $\Re s<1$, $b\not\in\MB{Z}$ follow analogously as in the proof for \tm{maintheorem}.
		The statement about meromorphicity and convergence on the \lhs follow directly from \lem{mellintransform:III}.
	}
	
	\section{Applications}\label{s:app}
	
	The immediate application of \tm{maintheorem} is the functional equation for the Hurwitz zeta function. In principal we have encountered it already in the proof of \tm{maintheorem}. We recall the polylogarithm $\on{Li}_s$ defined in \eq*{polylog}.
	
	\be[Hurwitz identity]{corollary}{c:Hurwitz}
	{
		Let $b\in\MB{R}$ and let $s\in\MB{C}\setminus\mbs{1}$. Then it holds
		\be{enumerate}{enumroman}
		{
			\item\label{c:Hurwitz:I}
			$
			{
				\zeta\mb{s,1-\mbs{b}}=\mf{\Gamma\mb{1-s}}{\mb{2\pi}^{1-s}}\mb{e^{\iunit \mf \pi 2 \mb{1-s}}\operatorname{Li}_{1-s}\mb{e^{\iunit 2\pi b}}+e^{-\iunit \mf \pi 2 \mb{1-s}}\operatorname{Li}_{1-s}\mb{e^{-\iunit 2\pi b}}}.
			}$
			\item\label{c:Hurwitz:II} $\mf{\sin\mb{\pi s}{\Gamma\mb{1-s}}\on{Li}_{1-s}\mb{e^{\iunit 2\pi b}}}{2^{-s}\pi^{1-s}}= e^{-\iunit\mf\pi 2s}\zeta\mb{s,1-\mbs{b}}+e^{\iunit\mf{\pi} 2 s}\zeta\mb{s,1-\mbs{-b}}.$
		}
	}
	\be*{proof}
	{
		Taking $\Re s<0$ and $a>0$ in \tm{maintheorem}, making use of the second identity, we find
		\be*{align*}
		{
			\zeta\mb{s,1-\mbs{b}}&=-\mf{\Gamma\mb{1-s}}{\mb{2\pi}^{1-s}\iunit}\mb{\mb{-\iunit }^s\on{Li}_{1-s}\mb{e^{\iunit 2\pi b}}-{\iunit }^s\on{Li}_{1-s}\mb{e^{-\iunit 2\pi b}}}\\
			&=-\mf{\Gamma\mb{1-s}}{\mb{2\pi}^{1-s}}\mb{e^{-\iunit\mf \pi 2 \mb{s+1}}\on{Li}_{1-s}\mb{e^{\iunit 2\pi b}}-e^{\iunit\mf \pi 2 \mb{s-1}}\on{Li}_{1-s}\mb{e^{-\iunit 2\pi b}}}\\
			&=\mf{\Gamma\mb{1-s}}{\mb{2\pi}^{1-s}}\mb{e^{\iunit\mf \pi 2 \mb{1-s}}\on{Li}_{1-s}\mb{e^{\iunit 2\pi b}}+e^{-\iunit\mf \pi 2 \mb{1-s}}\on{Li}_{1-s}\mb{e^{-\iunit 2\pi b}}}.
		} 
		To see the second equality, we compute
		\be*{align*}
		{
			\zeta\mb{s,1-\mbs{b}}-&e^{-\iunit\pi\mb{1-s}}\zeta\mb{s,1-\mbs{-b}}=\\
			\mf{\Gamma\mb{1-s}}{\mb{2\pi}^{1-s}}&\mb{e^{\iunit\mf \pi 2 \mb{1-s}}-e^{-\iunit \mf {3\pi} 2 \mb{1-s}}}\on{Li}_{1-s}\mb{e^{\iunit 2\pi b}}=\\
			&2\iunit e^{-\iunit{\pi} \mb{\mf 12-s}}\sin\mb{\pi\mb{1-s}}\mf{\Gamma\mb{1-s}}{\mb{2\pi}^{1-s}}e^{-\iunit\mf{\pi} 2s}\on{Li}_{1-s}\mb{e^{\iunit 2\pi b}}.
		}
		By meromorphic continuation for the Hurwitz zeta function, using \tm{maintheorem}, we conclude for all $s\in\MB{C}$.
	}
	
	\be[Functional equation for Riemann zeta function]{corollary}{c:RZ_FE}
	{
		Let $s\in\MB{C}\setminus\mbs{1}$. Then it holds
		\be*{align*}
		{
			\zeta\mb{s}=2^s \pi^{s-1}\Gamma\mb{1-s}\sin\mb{\mf\pi 2 s}\zeta\mb{1-s}.
		}
	}
	\be*{proof}
	{
		We take $b=0$ in \co[I]{Hurwitz}.
	}
	
	Since the convergence of the sum in $s$ such that $\Re s>0$ is not absolute if $b=0$ this separates the identity for the Riemann zeta function from the general identity for the Hurwitz zeta function (and therefore for every Dirichlet series with periodic coefficients). Hence we find an apparently non-trivial classification of zeros of the Riemann zeta function (on the critical strip):
	
	\be*{corollary}
	{
		Let $\rho\in\MB{C}$. Then for every $a>0$ it holds $\zeta\mb{\rho}=0$ if and only if
		\be*{align*}
		{
			&\mf 12\mb{\mb{d\mf {a\rho}{1-\rho}-\mf 12}d^{-\rho}+R_{a,0,d}\mb{\rho}}=\\
			&\sum_{k=1}^\infty \mf {ad^{1-\rho}\rho}{1-\rho}{}_1F_2\mb{\mf 12\mb{1-\rho};\mf 32;1+\mf 12\mb{1-\rho};-\pi^2a^2d^2k^2}+a^\rho\mf{\Gamma\mb{1-\rho}\sin\mb{\mf {\pi}2\rho}}{\mb{2\pi}^{1-\rho}k^{1-\rho}}.
		}
	}
	\be*{proof}
	{
		We take $b=0$ in \tm{maintheorem}.
	}
	
	\be*{corollary}
	{
		Let $\rho\in\MB{C}\setminus\MB{Z}$. Then for every $a\in\MB{R}\setminus\mbs{0}$ and for every $d>0$ it holds $\zeta\mb{\rho}=0$ if and only if
		\be*{align*}
		{
			&\mf{\log 2}{d^\rho}+1_{|\MB{R}\setminus\MB{Z}}\mb{ad}\mf {\log\mb{\sin\mb{\pi\mb{ad}}}^2}{2 d^\rho}+\mf {1}{\rho d^\rho}\\
			&+\mf {1}{ d^\rho}1_{|\MB{Z}}\mb{ad}\mb{\mf {1}{\rho}-{\psi\mb{1+\rho}-\gamma}+\mf 12\log\mb{\pi ad}^2}\\
			&+\mf{1}{d^{\rho}}\sum_{\substack{k\in\MB{Z}\\
					k\not\in\mbs{0,ad}}}\mf{ad}{k}\mf{{}_2F_1\mb{1,1-\rho;2-\rho;\mf{ad}{k}}}{{1-\rho}}+\iunit{\pi\mbm{a}^\rho}\zeta\mb{\rho,n}=\\			
			&-\sum_{k=1}^\infty\mf{1}{d^\rho k}{}_1F_2\mb{-\mf\rho 2;\mf 12,1-\mf\rho  2;-\pi^2a^2d^2k^2}-\mbm{a}^\rho\mf{\Gamma\mb{1-\rho}\cos\mb{\mf \pi 2 \rho}}{\mb{2\pi}^{-\rho}k^{1-\rho}},
		}
		and $n\in\MB{Z}$ is the smallest integer satisfying $n>\mbm{a}d$.
	}
	\be*{proof}
	{
		We take $b=0$ in \tm{maintheorem:B}.
	}
	
	\be{lemma}{l:recurrence:odd}
	{
		Let $n\in\MB{N}$, let $a\in\MB{R}\setminus\mbs{0}$, $b\in\MB{R}$ and let $d>0$. Then it holds
		
		\be*{align*}
		{
			\mf \pi n&\mb{ a^{-n}\zeta\mb{-n,1-\mbs{b}}-\mb{-\mf {adn}{1+n}+\mbs{b}-\mf 12}d^n+R_{a,b,d}\mb{-n}}=\\
			&d^n\sum_{m=1}^n\mf{\mb{1-n}_{m-1}}{2\mb{2\pi ad}^m}\mb{\mb{-\iunit}^{m+1}\on{Li}_{m+1}\mb{e^{\iunit 2\pi \mb{ad+b}}}+{\iunit}^{m+1}\on{Li}_{m+1}\mb{e^{-\iunit 2\pi \mb{ad+b}}}}.
		}
		In particular, if $ad\in\MB{Z}$ it holds
		\be*{align*}
		{
			\mf \pi n&\mb{\mb{\mf {n}{1+n}-\mbs{b}+\mf 12}d^n+R_{a,b,d}\mb{-n}}=\\
			&d^n\sum_{m=1}^{n-1}\mf{\mb{1-n}_{m-1}}{2\mb{2\pi ad}^m}\mb{\mb{-\iunit}^{m+1}\on{Li}_{m+1}\mb{e^{\iunit 2\pi b}}+{\iunit}^{m+1}\on{Li}_{m+1}\mb{e^{-\iunit 2\pi {b}}}}.
		}
	}
	
	\be*{proof}
	{
		From \tm{convergence} we know for every $s\in\MB{C}$, $a\in\MB{R}\setminus\mbs{0}$, $b,d\in\MB{R}$
		\be*{align*}
		{
			\widetilde{I}_s\mb{a,b,0,d}\coloneqq\sum_{k=1}^\infty F_{s,k}\mb{a,b,0,d}
		}
		converges. Now $-\widetilde{I}_{-n}\mb{a,b,0,d}$ evaluates to the \lhs of the claimed identity, using \tm{maintheorem}. By \lem{mellintransform:II:recurrence}, $F_{s,k}$ satisfies for every $n\in\MB{N}$
		\[F_{-n,k}\mb{a,b,0,d}=-d^n\sum_{m=0}^{n-1}\mf{\mb{1-n}_m}{2\mb{2\pi ad}^{m+1}}\mf{\mb{-\iunit}^{m+2}e^{\iunit 2\pi\mb{ad+b}k}-\iunit^{m+2}e^{-\iunit 2\pi\mb{ad+b}k}}{k^{m+2}},\]
		recalling holomorphicity of $F_{s,k}$ in $s$ and $\mb{1-n}_n=0$. Then the summation in $k$ yields the first claim. To see the statement in the event $ad\in\MB{Z}$, we note
		\be*{align*}
		{
			d^n\mf{\mb{1-n}_{n-1}}{2\mb{2\pi ad}^n}&\mb{\mb{-\iunit}^{n+1}\on{Li}_{n+1}\mb{e^{\iunit 2\pi \mb{ad+b}}}+{\iunit}^{n+1}\on{Li}_{n+1}\mb{e^{-\iunit 2\pi\mb{ad+b}}}}=\\
			&\mf{\Gamma\mb{n+1}}{2\mb{2\pi a}^nn}\mb{{\iunit}^{n+1}\on{Li}_{n+1}\mb{e^{\iunit 2\pi {b}}}+\mb{-\iunit}^{n+1}\on{Li}_{n+1}\mb{e^{-\iunit 2\pi {b}}}},
		}
		and conclude, combining this identity with \co{Hurwitz} and the first claim.
	}
	
	Let $n\in\MB{N}$. Set as the generalized log-sine integral for every $a\in\MB{R}\setminus\mbs{0}$, $b,d\in\MB{R}$
	\be*{align*}
	{
		\on{Ls}^{\mb{n-3}}_{n-1}\mb{a,b,d}\coloneqq-\mf 12\int_0^d \log\mb{4\mb{\sin\mb{\pi \mb{ax+b}}}^2}x^{n-1}\on{d}x.
	}
	\be{lemma}{l:recurrence:even}
	{
		For every $s\in\MB{C}$, $a\in\MB{R}\setminus\mbs{0}$, $b\in\MB{R}$, $d>0$
		\be*{align*}
		{
			&\on{Ls}^{\mb{n-3}}_{n-1}\mb{a,b,d}=\mf{\Gamma\mb{n}}{2a^n \mb{2\pi}^n}\mb{{\iunit}^{n}\on{Li}_{n+1}\mb{e^{\iunit 2\pi b}}+\mb{-\iunit}^{n}\on{Li}_{n+1}\mb{e^{-\iunit 2\pi b}}}\\
			&+d^n\sum_{m=1}^n\mf{\mb{1-n}_{m-1}}{2\mb{2\pi ad}^m}\mb{\mb{-\iunit}^m\on{Li}_{m+1}\mb{e^{\iunit 2\pi \mb{ad+b}}}+{\iunit}^m\on{Li}_{m+1}\mb{e^{-\iunit 2\pi \mb{ad+b}}}}.
		}
		In particular, if $ad\in\MB{Z}$ it holds
		\be*{align*}
		{
			\on{Ls}^{\mb{n-3}}_{n-1}&\mb{a,b,d}=d^n\sum_{m=1}^{n-1}\mf{\mb{1-n}_{m-1}}{2\mb{2\pi ad}^m}\mb{\mb{-\iunit}^m\on{Li}_{m+1}\mb{e^{\iunit 2\pi {b}}}+{\iunit}^m\on{Li}_{m+1}\mb{e^{-\iunit 2\pi {b}}}}.
		}
	}
	\be*{proof}
	{
		For every $a\in\MB{R}\setminus\mbs{0}$, $b\in\MB{R}$, $d>0$
		\be*{align*}
		{
			-\sum_{k=1}^\infty f_{s,k}\mb{a,b,\mf 14,d}=-\mf 12 \sum_{k=1}^\infty &\mb{f_{s,k}\mb{a,b,\mf 14,d}+\iunit f_{s,k}\mb{a,b,0,d}}\\
			&+\mb{f_{s,k}\mb{-a,-b,\mf 14,d}+\iunit f_{s,k}\mb{-a,-b,0,d}}
		}
		evaluates to the \lhs of the claimed identity if $s\in\MB{Z}_{<0}$, recalling \lem{mellintransform:I}. On the other hand, we find, combining \lem{mellintransform:I:recurrence}, taking advantage of \eq*{expidentity} in \lem{mellintransform:I}, with \lem{mellintransform:II:recurrence},
		\be*{align*}
		{
			g_{s,k}\mb{a,b,\mf 14}+\iunit g_{s,k}\mb{a,b,0}+&\mf{\mb{s+1}_n}{2\mb{\iunit2\pi ak}^n}\mb{f_{s+n,k}\mb{a,b,\mf 14,d}+\iunit f_{s+n,k}\mb{a,b,0,d}}\\
			=\mf{\mb{s+1}_n}{2\mb{\iunit2\pi ak}^n}&\mb{F_{s+n,k}\mb{a,b,\mf 14,d}+\iunit F_{s+n,k}\mb{a,b,0,d}}.
		}		
		For $s=-n$ the above \rhs vanishes by an analgous argument as in the proof for \lem{recurrence:odd}, 
		and therefore it holds
		\be*{align*}
		{
			\lim_{s\to -n}\mf{\mb{s+1}_n}{2\mb{\iunit2\pi ak}^n}&\mb{f_{s+n,k}\mb{a,b,\mf 14,d}+\iunit f_{s+n,k}\mb{a,b,0,d}}=\\
			&\quad\quad\quad\quad\quad\quad-g_{-n,k}\mb{a,b,\mf 14}-\iunit g_{-n,k}\mb{a,b,0}.
		}
		Evaluation of $-\sum_{k=1}^\infty f_{-n,k}\mb{a,b,\mf 14,d}$, using again the recurrence relation from \lem{mellintransform:I:recurrence} together with the identity of \eq*{expidentity} in \lem{mellintransform:I}, yields the first final claim. The second claim is analogous to the second claim from \lem{recurrence:odd}.
	}
	
	We end this section with the collection of hypergeometric identities which are  natural \wrt \se{result}. Notably, except for the extension from \lem{hyperidentity} of the identity in \eq*{hyperidentity} we did not make direct use of these but rather their incomplete form as a recurrence relation, in order to show meromorphic extension. Stated otherwise the following identities appear as the limiting cases of the recurrence relations stated in \se{result}.
	
	\be[Hypergeometric identities]{proposition}{p:hyperidentity1F1}
	{
		Let $s\in\MB{C}\setminus\MB{N}$ and let $a\in\MB{C}$. Then it holds
		
		\be{enumerate}{enumroman}
		{
			\item\label{p:hyperidentity1F1:I}${}_1F_1\mb{-s,1-s,\iunit 2\pi a}=e^{\iunit 2\pi a}{}_1F_1\mb{1,1-s,-\iunit 2\pi a}$ and,
			\item\label{p:hyperidentity1F1:II}recalling $f_{s,k}$ defined in \eq*{fk}, for every $k\in\MB{N}$, $b,c\in\MB{C}$, $d>0$ \be*{align*}
			{
				f_{s,k}\mb{a,b,c,d}=&\mf{d^{-s}}s\sin\mb{2\pi\mb{adk+bk+c}}{}_1F_2\mb{1;\mf 12-\mf s2,1-\mf s2;-\pi^2a^2d^2k^2}\\
				-\mf{2\pi ad^{1-s}k}{s\mb{1-s}}&\cos\mb{2\pi\mb{adk+bk+c}}{}_1F_2\mb{1;1-\mf s2,\mf 32-\mf s2;-\pi^2a^2d^2k^2},
			}
			\item\label{p:hyperidentity1F1:III}$\mf{a}{1-s}{}_2F_1\mb{1,1-s,2-s,-a}=\mf{\mb{1+a^{-1}}^{-1}}{1-s}a{}_2F_1\mb{1,1,1-s,\mb{1+a}^{-1}}.$
		}
	}
	\be*{proof}
	{
		From \lem{mellintransform:I:recurrence} we find
		\be*{align*}
		{
			&{}_1F_1\mb{-s,1-s,\iunit 2\pi a}-\mf{\mb{\iunit 2\pi a}^n}{\mb{s-n}_n}{}_1F_1\mb{-s+n,1-s+n,\iunit 2\pi a}\\
			&=e^{\iunit 2\pi a}\sum_{m=1}^n\mf{\mb{s-n}_m}{\mb{s-n}_n}\mb{\iunit 2\pi a}^{n-m}=e^{\iunit 2\pi a}\sum_{m=0}^{n-1}\mf{\mb{s-n}_{n-m}}{\mb{s-n}_n}\mb{\iunit 2\pi a}^m\\
			&=e^{\iunit 2\pi a}\sum_{m=0}^{n-1}\mb{1-s}_m^{-1}\mb{-\iunit 2\pi a}^m.
		}
		Hence taking the limit $\lim_{n\to \infty}$ on both sides, assuming $s$ is fixed away from the poles at $\MB{N}$ yields the first claim, noting $_1F_1\mb{-s+n,1-s+n,\iunit 2\pi a}$ is bounded for every $a$ and $n$ which follows from its representation as hypergeometric series. Claim \pp*[II]{hyperidentity1F1} follows by reversing the identity \eq*{expidentity} in \lem{mellintransform:I} together with claim \pp*[I]{hyperidentity1F1}. In particular, we first find
		\be*{align*}
		{
			2\iunit&\mb{f_{s,k}\mb{a,b,c,d}+f_{s,k}\mb{-a,-b,-c,d}}=\\
			&\mf{e^{\iunit 2\pi\mb{bk+c}}}{sd^sk}{}_1F_1\mb{-s;1-s;\iunit 2\pi adk}+\mf{e^{\iunit 2\pi\mb{-bk-c}}}{sd^sk}{}_1F_1\mb{-s;1-s;-\iunit 2\pi adk}.
		}
		Then by considering the corresponding hypergeometric series for ${}_1F_1$, utilizing
		\be*{align*}
		{
			\mf{\Gamma\mb{s+1+2m-2n}}{\Gamma\mb{s+1}}&=\mf{\mb{2-s}_{2n-2-2m}^{-1}}{\mb{s-1}s}=\mf{2^{m+1-n}}{\mb{s-1}s}\mb{1-\mf s2}_{n-1-m}^{-1}\mb{\mf 32-\mf s2}_{n-1-m}^{-1},\\
			\mf{\Gamma\mb{s+2+2m-2n}}{\Gamma\mb{s+1}}&=\mf{\mb{1-s}_{2n-2-2m}^{-1}}{s}=\mf{2^{m+1-n}}{s}\mb{\mf 12-\mf s2}_{n-1-m}^{-1}\mb{1-\mf s2}_{n-1-m}^{-1},
		}		
		the result follows. Claim \pp*[III]{hyperidentity1F1} follows, taking the limit $n\to\infty$ in the identity in \lem{mellintransform:III:recurrence} with $s$ away from the poles. Note that ${}_2F_1$ in there can grow at most exponentially, which follows from its integral representation.
	}
	
	\printbibliography[heading=none]
}